\documentclass[english,reqno]{amsart}
\usepackage[margin=1in]{geometry}

\usepackage{amsmath, appendix, ulem}
\usepackage[nobysame]{amsrefs}
\usepackage{amssymb, color}
\usepackage{amsfonts}
\usepackage{graphicx}
\usepackage{float}
\renewcommand{\tilde}{\widetilde}
\renewcommand{\hat}{\widehat}
\renewcommand{\bar}{\overline}
\newcommand{\eps}{\varepsilon}
\newcommand{\R}{\mathbb{R}}

\begin{document}
\title{Effects of internal dynamics on chemotactic aggregation of bacteria}
\date{\today}

\author{Shugo YASUDA}
\address{Graduate School of Information Science, University of Hyogo, 650-0047 Kobe, Japan}
\email{yasuda@gsis.u-hyogo.ac.jp}
\thanks{This work was supported by the Japan-France Integrated Action Program (SAKURA), Grant number JPJSBP120193219.}

\begin{abstract}
The effects of internal adaptation dynamics on the self-organized aggregation of chemotactic bacteria are investigated by Monte Carlo (MC) simulations based on a two-stream kinetic transport equation coupled with a reaction-diffusion equation of the chemoattractant that bacteria produce.

A remarkable finding is a nonmonotonic behavior of the peak aggregation density with respect to the adaptation time; more specifically, aggregation is the most enhanced when the adaptation time is comparable to or moderately larger than the mean run time of bacteria.
Another curious observation is the formation of a trapezoidal aggregation profile occurring at a very large adaptation time, where the biased motion of individual cells is rather hindered at the plateau regimes due to the boundedness of the tumbling frequency modulation.

Asymptotic analysis of the kinetic transport system is also carried out, and a novel asymptotic equation is obtained at the large adaptation-time regime while the Keller-Segel type equations are obtained when the adaptation time is moderate.
Numerical comparison of the asymptotic equations with MC results clarifies that trapezoidal aggregation is well described by the novel asymptotic equation, and the nonmonotonic behavior of the peak aggregation density is interpreted as the transient of the asymptotic solutions between different adaptation time regimes.
\end{abstract}

\maketitle

\section{Introduction}
The collective motion of chemotactic bacteria, such as \textit{Escherichia coli}, stems from, at the individual level, continuous reorientations by runs and tumbles.
It has been established that the length of a run is determined by a stiff response to the temporal variation of extracellular chemical cues via an intracellular signal transduction pathway.
The chemotactic response and the intracellular signal transduction pathway for \textit{E. coli} have been extensively studied by various authors, and sophisticated mathematical models have been proposed.~\cites{BL1997,KLBTS2005,TSB2008,JOT2010,H2012}
However, the multiscale mechanism between intracellular signal transduction, individual chemotactic motion, and collective dynamics of cells is not yet well understood.
Currently, engineered bacteria (or genetically modified bacteria) are utilized in a variety of industrial fields involving, for example, food, agriculture, medicine, and the environment.
Understanding the multiscale mechanism can contribute to further advances in industrial technology to control the collective motion of cells.

Kinetic transport models have been proposed to describe the multiscale mechanism in the collective motion of cells;
a kinetic transport model describing the velocity jump process in the run-and-tumble motion of \textit{E. Coli} was first proposed in Ref.~\cites{ODA1998}, and it was then further developed to involve more detailed chemosensory systems~\cites{DS2005,EO2004,SWOT2012}.
Although the chemosensory system involves complicated biochemical reaction networks, in the simplified description, it can be constituted by two essential steps, i.e., a rapid response to an external signal change called ``excitation'' and a subsequent slow ``adaptation'', in which the internal state returns to the baseline, allowing the cell to respond to a further external signal change.~\cite{SPO1997}
A kinetic transport equation involving the simplified description of the excitation and adaptation dynamics has been proposed in Refs.~\cite{EO2004} and \cite{EO2007}, where two internal state variables are introduced in the model.
Since the excitation dynamics are much faster than the adaptation dynamics, one can integrate the fast variable related to the excitation and derive the kinetic transport equation involving only a single internal state variable related to the slow adaptation dynamics~\cite{SWOT2012}, i.e.,
\begin{equation}\label{eq_p0}
\partial_t p+\mathbf{v}\cdot\nabla_x p + \partial_m[F(m,S) p]=\mathcal{Q}[m,S](p),
\end{equation}
where $p(t,\mathbf{x},\mathbf{v},m)$ is the density of cells with velocity $\mathbf{v}\in \mathbb{V}$ and internal state $m>0$ at time $t>0$ and position $\mathbf{x} \in \R^d$.
Here, on the left-hand side, the $x$-divergence term describes the change in density due to the ``run'' of the bacteria, and the $m$-derivative term describes the evolution of the internal state $m$ at the rate of change $F(m,S)$, where $S(t,x)$ is the concentration of the extracellular chemical cue. On the right-hand side, $\mathcal{Q}[m,S](p)$ is the tumbling operator described as
\begin{equation}\label{eq_RHS}
\mathcal{Q}[m,S](p)=\frac{1}{||\mathbb{V}||}\int_\mathbb{V}
[\lambda(m,S,{\bf v},{\bf v}')p(t,{\bf x},{\bf v}',m)-\lambda(m,S,{\bf v}',{\bf v})p(t,{\bf x},{\bf v},m)]d{\bf v}',
\end{equation}
where $\lambda(m,S,{\bf v},{\bf v}')$ denotes the tumbling frequency describing the reorientation from velocity ${\bf v}'$ to the new velocity ${\bf v}$.
The velocity space $\mathbb{V}$ is the bounded domain of $\R^d$ and $||\mathbb{V}||=\int_V\,dv$.

Since the bacteria communicate with each other via the extracellular chemical cues they produce, to describe the observed self-organization phenomena occurring in a population of chemotactic bacteria, for example, in Refs. ~\cites{BB1991,WTMMBB1995,MBBO2003}, chemoattractant equations must be coupled with the kinetic transport model (\ref{eq_p0}).
In this study, we consider a single species of chemical cues whose concentration $S(t,\mathbf{x})$ is described as
\begin{equation}\label{eq_S}
\partial_t S = D_S\Delta S - a S + b\rho,
\end{equation}
where $D_S$ is the diffusion coefficient of the chemical cue, $a$ is the degradation rate of the chemical cue, $b$ is the production rate of the chemical cue by bacteria, and $\rho(t,{\bf x})=\int_\mathbb{V}\int_0^\infty p(t,{\bf x},{\bf v},m)dmd{\bf v}$ is the population density of bacteria.

In the kinetic transport model, the microscopic characteristics at the individual level are involved in the tumbling frequency $\lambda(m,S,{\bf v},{\bf v}')$ and the rate of change of the internal state $F(m,S)$.
Thus, by specifying the mathematical formulas for $\lambda(m,S,{\bf v},{\bf v}')$ and $F(m,S)$, one can address the multiscale mechanism between the intracellular adaptation dynamics, individual chemotactic motion, and collective dynamics of cells in the self-organization phenomena.
One can also derive macroscopic models for the population density of bacteria, e.g., Keller-Segel (KS)-type systems~\cites{KS1970,KS1971,HP2009}, and kinetic transport equations without internal state variables~\cite{DS2005} by using moment closure or asymptotic analysis of Eqs. (\ref{eq_p0})--(\ref{eq_S}).~\cites{PTV2016,PVW2018,ST2017,X2015}

Investigations of the aggregation of chemotactic bacteria based on kinetic transport models have been carried out in various studies.
For example, in Ref.~\cite{SM2011}, the aggregation of chemotactic bacteria under a given concentration gradient of a chemical cue was investigated based on the kinetic transport model with internal states, and the volcano-like (bimodal) aggregation of \textit{E. coli} observed in an experiment~\cite{MBBO2003} was numerically reproduced in one-dimensional space.
In Ref.~\cite{XXT2018}, the concentric stripe patterns formed by engineered \textit{E. coli}~\cite{Liuetal2011} were reproduced numerically, and the role of intracellular signal transduction in stripe pattern formation was clarified.
In Ref.~\cite{PY2018}, the instability of the kinetic transport model describing colony pattern formation over a long period of time due to proliferation was investigated, and stiff-response-induced instability was uncovered at the kinetic level.
Additionally, in Ref.~\cite{R2020}, the role of the hydrodynamic interaction in the self-organized aggregations was numerically investigated by using a Monte Carlo method related to a kinetic transport model without internal states.
These studies have established that kinetic transport models are useful for elucidating the multiscale mechanism in the collective motion of chemotactic bacteria.
However, the multiscale mechanism between collective motions and internal state dynamics in self-organized aggregation has yet to be clarified.

In this paper, we investigate the self-organized aggregation of chemotactic bacteria in one-dimensional space based on a two-stream kinetic transport model with an internal state.
In contrast to the previous study~\cite{PY2018}, this paper concerns the internal dynamics of chemotactic bacteria and considers self-organized aggregation, which may occur in a rather short period of time without proliferation.
In particular, we focus on the effect of the adaptation time on the instability and aggregation behavior.

In the following text, the problem and the basic equations are given in Sec.~\ref{sec.problem}.
In Sec. ~\ref{sec.numeric}, numerical analyses are carried out for a wide range of adaptation times by using a Monte Carlo (MC) method, which is an extension of the MC method previously developed in Refs.~\cites{Y2017,VY2020}.
In Sec.~\ref{sec.asymptotic}, we formally carry out asymptotic analysis of the kinetic transport model at different scalings of the adaptation time and derive a KS-type model and a novel asymptotic equation involving the internal state variable.
The asymptotic behavior is also numerically investigated over a wide range of adaptation times, through which a suitable parameter regime for the KS-type system and a remarkable numerical solution in the novel asymptotic regime are uncovered.
Finally, a summary and perspectives are given in Sec.~\ref{sec.summary}.

\section{Problem and formulation}\label{sec.problem}
We consider the chemotactic bacteria moving in positive and negative directions with a constant speed $V_0$, i.e., $v=\{-V_0,V_0\}$, in one-dimensional space $x\in[0,L]$ with periodic boundary conditions.
Initially, the bacteria are uniformly distributed, and the internal state $m$ is in the equilibrium state at $m=M(S)$, where $M(S)$ denotes the equilibrium internal state determined by the extracellular chemical concentration.
The chemical concentration $S(t,x)$ is also uniformly distributed in the initial state.

For the internal state dynamics, we consider the following linear adaptation model:
\begin{equation}\label{eq_m}
\frac{dm}{dt}=F(m,S)=\frac{M(S)-m}{\tau},
\end{equation}
where $\tau>0$ denotes the characteristic adaptation time.
We also assume that the bacteria tumble (i.e., change in moving direction) depends only on the deviation of the internal state $m$ from the equilibrium state $M(S)$, $M(S)-m$,~\cite{SWOT2012}:
\begin{equation}
\lambda(m,S,v,v')=\lambda_0\Lambda\left(\frac{M(S)-m}{\delta}\right),
\end{equation}
where $\lambda_0>0$ is the mean tumbling frequency, $\Lambda(\frac{M(S)-m}{\delta})>0$ denotes the modulation of the tumbling frequency, and $\delta>0$ denotes the stiffness of the chemotactic response.
In this study, we consider the following modulation function:
\begin{equation}\label{eq_Lambda}
\Lambda_\delta(y)=\Lambda\left(\frac{y}{\delta}\right),\quad
\Lambda(y)=1-R(y),\quad R(y)=\frac{\chi y}{\sqrt{1+y^2}},
\end{equation}
where $0<\chi<1$ denotes the modulation amplitude and $\delta$ denotes the stiffness of the chemotactic response.

Then, the density of bacteria with positive and negative velocities, $p^\pm(t,x,m)$, is described by the following two-stream kinetic transport equation with the internal state:
\begin{equation}\label{eq_ppm}
\partial_t p^\pm \pm V_0\partial_x p^\pm+
\partial_m\left(
\frac{M(S)-m}{\tau} p^\pm
\right)	
=\pm\frac{\lambda_0}{2}\Lambda_\delta\left(M(S)-m\right)(p^--p^+).
\end{equation}
By introducing the nondimensional variables
$$
\hat x=x/L_0,\quad \hat t=t/t_0,\quad \hat v=v/V_0,
$$
where $L_0$, $V_0$, and $t_0$ are the characteristic length, speed, and time, respectively, Eq. (\ref{eq_ppm}) is written in nondimensional form as
\begin{equation}\label{eq_p}
\hat \sigma \partial_{\hat t} \hat p^\pm \pm \partial_{\hat x} \hat p^\pm+
\partial_m\left(
\frac{M(S)-m}{\hat \tau} p^\pm
\right)	
=\frac{\hat \lambda_0}{2}\Lambda_\delta\left(M(S)-m\right)(\hat p^\mp-\hat p^\pm).
\end{equation}
Here, the nondimensional parameters $\hat \sigma$, $\hat \lambda_0$, and $\hat \tau$ are defined as
\begin{equation}\label{def_parap}
\hat \sigma=L_0/(t_0V_0),\quad \hat \lambda_0=\lambda_0/(V_0/L_0),\quad\hat \tau=\tau/(L_0/V_0).
\end{equation}
We also define $\hat p^\pm=p^\pm/\rho_0$, where $\rho_0$ is the initial population density of bacteria.
The characteristic time $t_0$ can be chosen arbitrarily depending on the time scale with which we address the problem.
For example, a typical choice of $t_0$ is $t_0=L_0/V_0$, which gives $\sigma$=1 and makes the kinetic transport equation (\ref{eq_p}) simpler, reducing one free parameter.
Another typical choice of $t_0$ is $t_0=(\lambda_0 L_0^2)/V_0^2$, which denotes the characteristic diffusion time and enables us to derive macroscopic continuum-limit equations such as KS-type models by asymptotic analysis of the kinetic transport equation.
In this paper, we address both time scales; i.e., we set $t_0=L_0/V_0$ for MC simulations of the kinetic transport model, while we use the diffusion time scale $t_0=(\lambda_0 L_0^2)/V_0^2$ for the asymptotic analysis to investigate the behaviors in the near-continuum regime.

By the same token, the nondimensional form of (\ref{eq_S}) is written as
\begin{equation}\label{eq_St}
\hat \sigma_S \partial_{\hat t} \hat S = \hat D_S \partial_{\hat x\hat x}\hat S- \hat S + \hat \rho,
\end{equation}
where the nondimensional quantities are defined as
\begin{equation}\label{def_paraS}
\hat \sigma_S=1/(at_0),\quad \hat D_S=D_S/(aL_0^2),\quad \hat S=S/(b\rho_0/a),\quad \hat \rho=\rho/\rho_0,
\end{equation}
and the population density $\hat \rho$ is calculated as
\begin{equation}\label{eq_rho}
\hat \rho(x,t)=\int_{0}^\infty\frac{\hat p^+(t,x,m)+\hat p^-(t,x,m)}{2}dm.
\end{equation}

In this study, we fix the nondimensional diffusion constant $\hat D_S=1$ (although the notation $D_S$ retains in the following equations for generality).
This indicates that the characteristic length $L_0$ denotes the diffusion length of chemoattractant $S$ within the degradation time $a^{-1}$, i.e., 
\begin{equation}\label{eq_L0}
L_0=\sqrt{\frac{D_S}{a}}.
\end{equation}

In the rest of the paper, unless otherwise stated, all quantities are written in nondimensional forms, and we drop the hat signs on the variables and parameters for simplicity.

It is convenient to introduce the new internal state variable $y=M(S)-m \in \mathbb{R}$ and change the variable to $f^\pm(t,x,y=M(S)-m)=p^\pm(t,x,m)$.
Then, we obtain
\begin{equation}\label{eq_f}
\sigma\partial_t f^\pm\pm\partial_xf^\pm
+\partial_y\left\{
\left(D^\pm_t M(S)-\frac{y}{\tau}
\right)f^\pm\right\}
=\pm\frac{\lambda_0\Lambda_\delta(y)}{2}(f^--f^+),
\end{equation}
where $D^\pm_t$ denotes the material derivative defined as $D_t^\pm=\sigma\partial_t\pm\partial_x$.
Here, the $y$-derivative term of (\ref{eq_f}) describes the change in the internal state variable $y$, where $D_t^\pm M(S)$ denotes the temporal variation of the extracellular chemical cue sensed by bacteria moving in positive and negative directions.
Since \textit{E. coli} cells respond to the spatial gradient of the logarithmic extracellular chemical concentration~\cites{BB1974,KJTW2009}, we model the logarithmic sensing by
\begin{equation}\label{def_DtM}
D_t^\pm M(S)= D_t^\pm \ln S(x,t)=\frac{\sigma\partial_t S\pm\partial_xS}{S}.
\end{equation}
Thus, the change in the internal state $y$ of each bacterium is described as
\begin{equation}\label{eq_doty}
\dot y=\frac{D_t^\pm S}{S}-\frac{y}{\tau},
\end{equation}
where $D_t^\pm S$ denotes the temporal variation of extracellular chemical cues sensed by bacteria along their moving trajectory.

It is clearly seen that the uniform state $f^\pm=\delta(y=0)$ and $S=\rho=1$, where $\delta(y)$ is the Dirac delta function, solves the system of Eqs.~(\ref{eq_St}) and (\ref{eq_f}).
We investigate the instability of the uniform state by the Monte Carlo code explained in Sec. \ref{sec.method}.

\section{Numerical analysis}\label{sec.numeric}

\subsection{Monte Carlo method}\label{sec.method}
The one-dimensional space $0\le x \le L$ is divided into the uniform mesh system $x_i=\Delta x\times i$ ($i=0,\cdots,I$) with the mesh width $\Delta x=L/I$, where $I$ is the number of mesh intervals.
Initially, 
Monte Carlo (MC)
 particles are uniformly distributed in each mesh interval with the equilibrium internal state at $y=0$.
The velocities of each MC particle, $v=\pm 1$, are randomly determined.
The chemical concentrations in each mesh interval $x\in[x_i,x_{i+1}]$ ($i=0,\cdots,I-1$), $S_i$, are also uniformly given at the initial state, i.e., $S_i^0=1$.

Then, the position $r^k_l$, velocity $v_l^k$, and internal state $y^k_l$ of the $l$th MC particle at time $t=k\Delta t$ are determined as follows:
\begin{enumerate}
\item Each MC particle moves as
\begin{equation}
r^k_l=r^{k-1}_l+v_l^{k-1}\Delta t.
\end{equation}
\item Population density in the $i$th mesh interval $x\in [x_i,x_{i+1}]$, $\rho^k_i$, is calculated as
\begin{equation}
\rho_i^k=\frac{1}{\bar N} \sum_{l=0}^N\int_{x_i}^{x_{i+1}}\delta(x-r_l^k)dx,
\end{equation}
where $\bar N$ is the number of MC particles in each mesh interval in the uniform state. Thus, the total number of MC particles is given by $I\times \bar N$.
\item Concentration of chemical cues in the $i$th mesh interval, $S^k_i$, is calculated explicitly as
\begin{equation}\label{eqS_disc}
\sigma_S\frac{S_i^{k}-S_i^{k-1}}{\Delta t}=\frac{D_S}{\Delta x^2} (S^{k-1}_{i+1}-2S^{k-1}_i+S^{k-1}_{i-1})-S_i^{k-1}+\rho_i^k.
\end{equation}
At the boundaries $x$=0 and $L$, we consider the periodic condition, i.e., $S_{-1}=S_{I-1}$ and $S_I=S_0$.

\item Internal state of the $l$th MC particle, $y^k_l$, is updated by following Eq.~(\ref{eq_doty}) as
\begin{equation}\label{eq_yupdate}
\frac{y^k_l-y^{k-1}_l}{\Delta t}=
\frac{S^k_{(l)}-S^{k-1}_{(l)}}{\Delta t {S^{k-1}_{(l)}}}-\frac{y^k_l}{\tau},
\end{equation}
where $S^k_{(l)}$ denotes the local concentration of chemical cues at the position of the $l$th MC particle $x=r^k_l$, i.e., $S^k_{(l)}=S(k\Delta t,r_l^k)$, and is calculated by linear interpolation:
\begin{equation}
S^k(l)=\left\{
\begin{array}{cc}
S_i^k+\frac{S^k_i-S^k_{i-1}}{\Delta x}(r_l^k-x_i-\frac{\Delta x}{2}),&\mathrm{if}\quad x_i\le r_l^k< x_i+\frac{\Delta x}{2},\\
S_i^k+\frac{S^k_{i+1}-S^k_{i}}{\Delta x}(r_l^k-x_i-\frac{\Delta x}{2}),&\mathrm{if}\quad x_i+\frac{\Delta x}{2}\le r_l^k < x_{i+1}.
\end{array}
\right.
\end{equation}
We note that in Eq.~(\ref{eq_yupdate}), the pathway derivative $D^\pm_t S$ in Eq.~(\ref{eq_doty}) is given by the rate of change of $S$ sensed by each bacterium, i.e., $(S_{(l)}^k-S_{(l)}^{k-1})/\Delta t$.
	\item Tumbling of the $l$th MC particle is decided by the probability 	$\frac{\Delta t\lambda_0}{2}\Lambda_\delta(y_l^k)$.
	\item The particles that decide to make tumbles change their velocities as $v_l^{k-1}\rightarrow v_l^k=-v_l^{k-1}$, and other particles retain their velocities.
\end{enumerate}

This MC method was applied for aggregation under a given constant spatial gradient of chemical cues, i.e., $\partial_x M(S)=$const., in Ref.~\cite{VY2020}, and the accuracy of the MC method was confirmed throughout the comparison to the asymptotic preserving schemes developed in the paper.

\subsection{Numerical results}\label{sec.numerics.result}
MC simulations are performed for various values of the mean tumbling frequency $\lambda_0$, the adaptation time $\tau$, and the stiffness of the chemotactic response $\delta$, while the diffusion constant $D_S=1$, the length of the periodic interval $L=10$, the modulation amplitude $\chi=0.5$, and the time scale parameters $\sigma$=$\sigma_S$=1 are fixed.
As is mentioned in Sec.~\ref{sec.problem}, since the diffusion constant is fixed as $D_S=1$, the length scale of the system corresponds to the diffusion length of chemoattractant defined as Eq.~(\ref{eq_L0}).

The number of mesh intervals $I$ and the average number of MC particles in each mesh interval $\bar{N}$ are set to $I$=50 and $\bar{N}$=28,800 except for the cases for $\lambda_0$=500 and 1000 at $\delta=0.01$, where $I=100$ and $\bar N=7,400$ are used.
The time step size is set to $\Delta t=10^{-3}$ for $\lambda_0<100$, $\Delta t=2\times 10^{-4}$ for $\lambda_0$=100 and 200, and $\Delta t=5\times 10^{-5}$ for $\lambda_0=$500 and 1000.

In the following, we introduce a new parameter $\alpha$, which is defined by the ratio of the adaptation time $\tau$ to the mean run time $\lambda_0^{-1}$, i.e., $\alpha=\lambda_0\tau$, and call it the relative adaptation time.

We also use the notation $t_\lambda$ for the scaled time $t_\lambda=t/(\lambda_0L^2)$.
The macroscopic population density $\rho$ in the following numerical results is time-averaged over the interval $\delta t_\lambda=0.05$ to remove the fluctuations caused by the Monte Carlo method.

\subsubsection{Instability and aggregation profile}
\begin{figure}[htbp]
\centering
\includegraphics[width=0.9\textwidth]{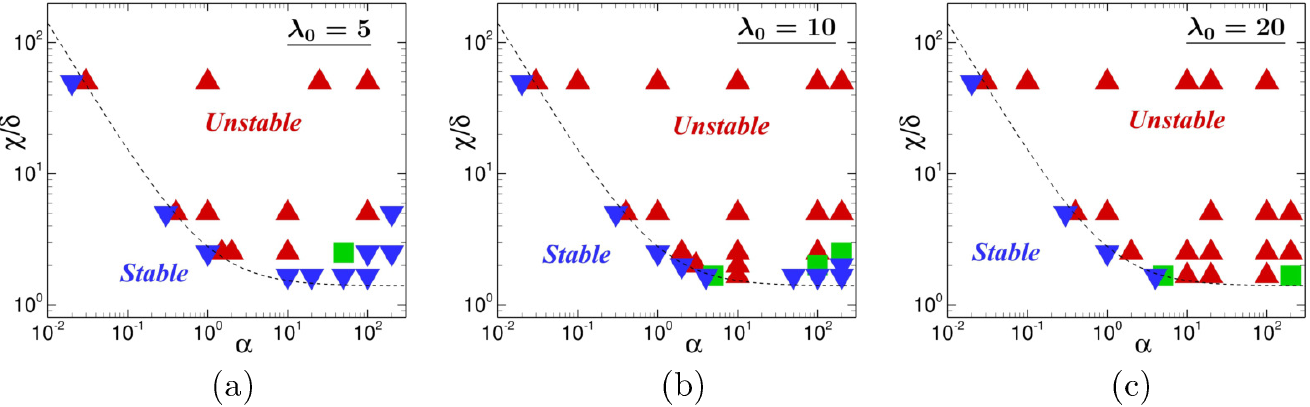}
\caption{
The instability diagrams with respect to the scaled adaptation time $\alpha$ and the stiffness $\chi/\delta$ at $\lambda_0=5$ (a), $\lambda_0=10$ (b) and $\lambda_0=20$ (c).
The upward triangles $\bigtriangleup$ show the results where stationary patterns are clearly observed, while the downward triangles $\bigtriangledown$ show the results where no evident patterns are observed (where the time average of the maximum deviation $\overline{\delta \rho}$ defined by (\ref{eq:deviation}) is less than 0.01).
See also Fig. \ref{fig.maxrho}.
The squares $\Box$ show the intermediate results where nonstationary sinusoidal waves with small amplitudes, i.e., $0.01<\overline{\delta \rho} <0.1$, are observed.
The dotted line shows the linear stability condition of the KS system, which is obtained in Sec. \ref{sec.ksinstability}.
Under the critical line, the uniform solution to the KS system is linearly stable.
}\label{fig.instability}
\end{figure}

Figure \ref{fig.instability} shows the instability diagrams with respect to the relative adaptation time $\alpha$ and the stiffness of the chemotactic response $\chi/\delta$ at different values of the tumbling frequency, i.e., $\lambda_0=5$ in (a), $\lambda_0=10$ in (b), and $\lambda_0=20$ in (c).
To confirm the stable state, we carried out a long-term simulation over $t_\lambda\in [0,T_\lambda]$ with $T_\lambda=10$ and measured the time average of the maximum deviation of the population density from the uniform state defined as
\begin{equation}\label{eq:deviation}
\overline{\delta \rho}=\frac{2}{T_\lambda}\int_{\frac{T_\lambda}{2}}^{T_\lambda}{\displaystyle\max_{x}}|\rho-1|dt_\lambda.
\end{equation}
The stable uniform states shown by the downward triangles $\bigtriangledown$ are confirmed when $\overline{\delta \rho}<0.01$, and the intermediate states shown by the squares $\Box$ are confirmed when $0.01<\overline{\delta \rho}<0.1$.

It is clear that when the relative adaptation time $\alpha$ is fixed, instability occurs when the stiffness $\chi/\delta$ is sufficiently large.
At small relative adaptation times, e.g., $\alpha \lesssim 1$, instability always occurs when the stiffness $\chi/\delta$ is larger than the critical value of the KS instability, which is obtained by Eq.~(\ref{eq_ksinstability}), and the transition between the stable and unstable regimes is very sharp.

It is also seen that the critical behavior for instability is not as affected by the mean tumbling frequency $\lambda_0$ in the small $\alpha$ regime.
However, in the large $\alpha$ regime, e.g., $\alpha > 10$, the instability behavior is significantly affected by the mean tumbling frequency $\lambda_0$, and the instability condition of the KS system (the dotted line in Fig.~\ref{fig.instability}) is no longer consistent with the MC results, especially at $\lambda_0=5$ (Fig.~\ref{fig.instability}(a)).

\begin{figure}[htbp]
\centering
\includegraphics[width=0.9\textwidth]{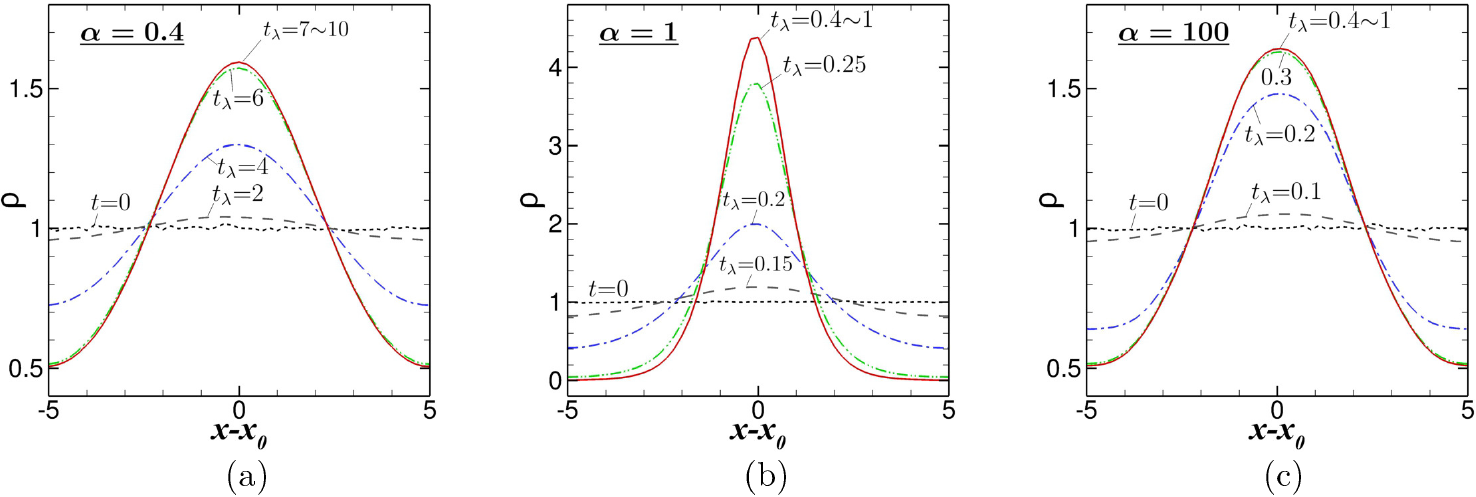}
\caption{
Time evolutions of the population density of bacteria with moderate stiffness $\delta=0.1$ at different values of the relative adaptation time, i.e., $\alpha=0.4$ (a), $\alpha=1$ (b), and $\alpha=100$ (c), at the mean tumbling frequency $\lambda_0=10$.
Here, $x_0$ represents the position where the chemical cue $S$ takes the maximum value in the stationary state.
We note that only the result at $t=0$ is calculated from the snapshot of the initial distribution of MC particles so that it involves relatively large fluctuations.
}\label{fig.macro.delta01}
\end{figure}
\begin{figure}[htbp]
\centering
\includegraphics[width=0.9\textwidth]{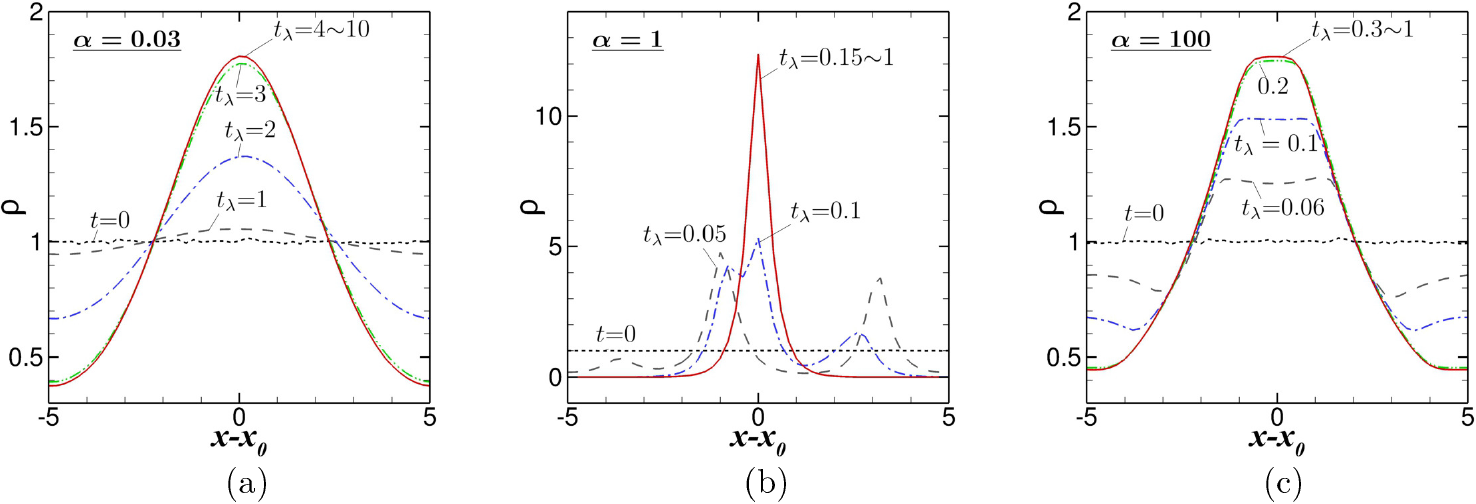}
\caption{
Time evolutions of the population density of bacteria with large stiffness $\delta=0.01$ at different values of the relative adaptation time, i.e., $\alpha=0.03$ (a), $\alpha=1$ (b), and $\alpha=100$ (c), at the mean tumbling frequency $\lambda_0=10$.
See also the caption in Fig. \ref{fig.macro.delta01}.
}\label{fig.macro.delta001}
\end{figure}

The aggregation profiles and their time evolution are shown in Figs.~\ref{fig.macro.delta01} and \ref{fig.macro.delta001}.
Figure \ref{fig.macro.delta01} shows the result at moderate stiffness $\delta=0.1$, while Figure \ref{fig.macro.delta001} shows the result at large stiffness $\delta=0.01$.
In both figures, the tumbling frequency $\lambda_0=10$ is fixed.
Initially, the population density $\rho$ is uniformly distributed with small fluctuations, whose amplitudes are at most 0.015.

The aggregation profiles are highly affected by the adaptation time and the stiffness.
At moderate stiffness $\delta=0.1$ (Fig.~\ref{fig.macro.delta01}), sinusoidal-like curves are generated in the stationary states at small and large adaptation times, i.e., $\alpha=0.4$ and 100, while at $\alpha=1$, a sharp aggregation profile is generated in the stationary state.

At large stiffness $\delta=0.01$ (i.e., Fig.~\ref{fig.macro.delta001}), very sharp aggregation, i.e., the spike-like aggregation profile, occurs at $\alpha=1$ in the stationary state.
On the other hand, interestingly, at large adaptation time $\alpha=100$ (see Fig. \ref{fig.macro.delta001}(c)), aggregation is not enhanced but is rather hindered at the central region of the aggregation profile such that the trapezoidal profile, where the plateau regimes appear at the top and bottom of the aggregate, is formed.
Later, we will see that this remarkable profile is obtained at a large stiffness when the adaptation time is as long as $\tau=O(\lambda_0)$.

\begin{figure}[htbp]
\centering
\includegraphics[width=0.8\textwidth]{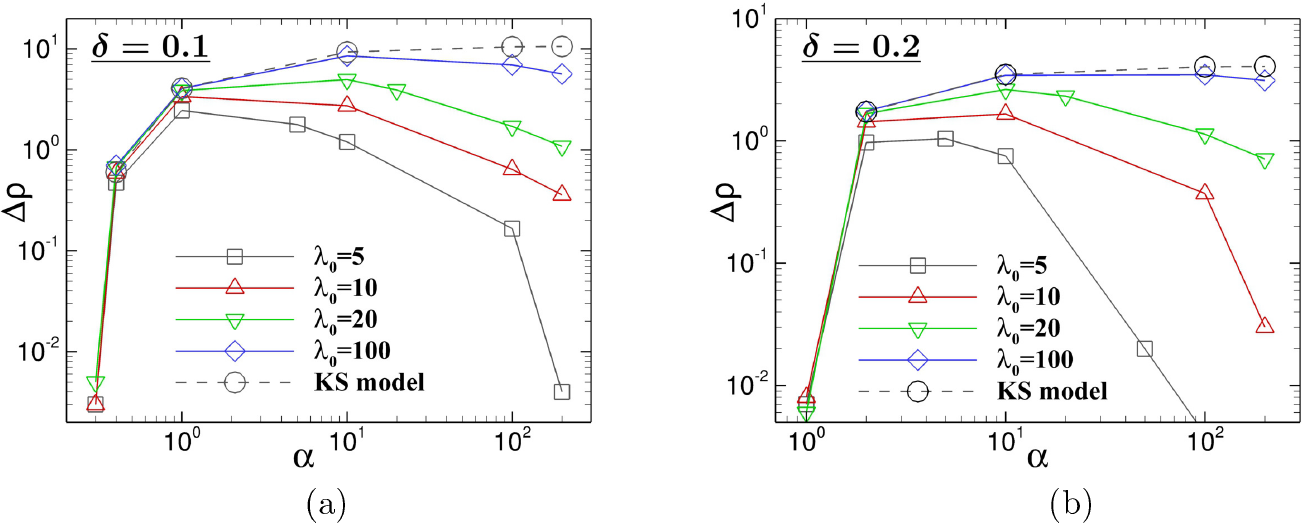}
\caption{
The dependency of the maximum aggregation density on the relative adaptation time $\alpha$ at different values of the mean tumbling frequency $\lambda_0$.
Figure (a) shows the result at $\delta=0.1$, and Figure (b) shows the result at $\delta=0.2$.
Here, the vertical axis shows the difference between the maximum aggregation density $\rho_\mathrm{max}$ and the uniform state $\rho=1$, i.e., $\Delta \rho=\rho_\mathrm{max}-1$.
The results of the KS system are obtained by numerical computation with the finite difference scheme on the staggered grid given in Ref.~\cite{CPY2018}.
}\label{fig.maxrho}
\end{figure}

Note that in both Figs.~\ref{fig.macro.delta01} and \ref{fig.macro.delta001}, the maximum aggregation density is not monotonically dependent on the adaptation time; aggregation is enhanced at moderate relative adaptation time $\alpha=1$.
The nonmonotonic dependency of the maximum aggregation density on the adaptation time is discussed in Fig.~\ref{fig.maxrho}.

The spatial profiles of chemical cue $S$ are shown in Fig. \ref{fig.xi}(c).
In contrast to the population density of bacteria $\rho$, the spatial distribution of $S$ is moderate and not significantly affected by either the relative adaptation time $\alpha$ or the stiffness $\delta$.

Note that the parameter sets used in Fig. ~\ref{fig.macro.delta01}(a) and Fig.~\ref{fig.macro.delta001}(a) are very close to but slightly above the linear stability condition of the KS system (the dotted line in Fig.~\ref{fig.instability}), where, although the time evolution is much slower than other cases, distinct sinusoidal aggregation occurs.
The sharp transition between the stable and unstable modes in the small relative adaptation-time regime (i.e., $\alpha<1$) is also observed in Fig.~\ref{fig.maxrho}.

The nonmonotonic behavior of the maximum aggregation density with respect to the relative adaptation time is seen in Figure \ref{fig.maxrho}.
It is clear that the transition from the stable to unstable modes is very sharp in the small $\alpha$ regime, i.e., $\alpha\lesssim 1$.
On the other hand, in the large $\alpha$ regime, the maximum aggregation density gradually decreases as $\alpha$ increases, and the slope of the decrease increases as $\lambda_0$ decreases.
This behavior is completely different from that in the KS system, where the maximum aggregation density monotonically increases and saturates to a certain value as the relative adaptation time $\alpha$ increases.

Note that at any fixed $\lambda_0$, there exists the optimal adaptation time to enhance aggregation in the regime $1<\alpha < \lambda_0$.
Importantly, this nonmonotonic behavior is a distinguished result obtained by the kinetic system but not by the KS system.

\subsubsection{Distribution of the internal state}

\begin{figure}[htbp]
\centering
\includegraphics[width=0.85\textwidth]{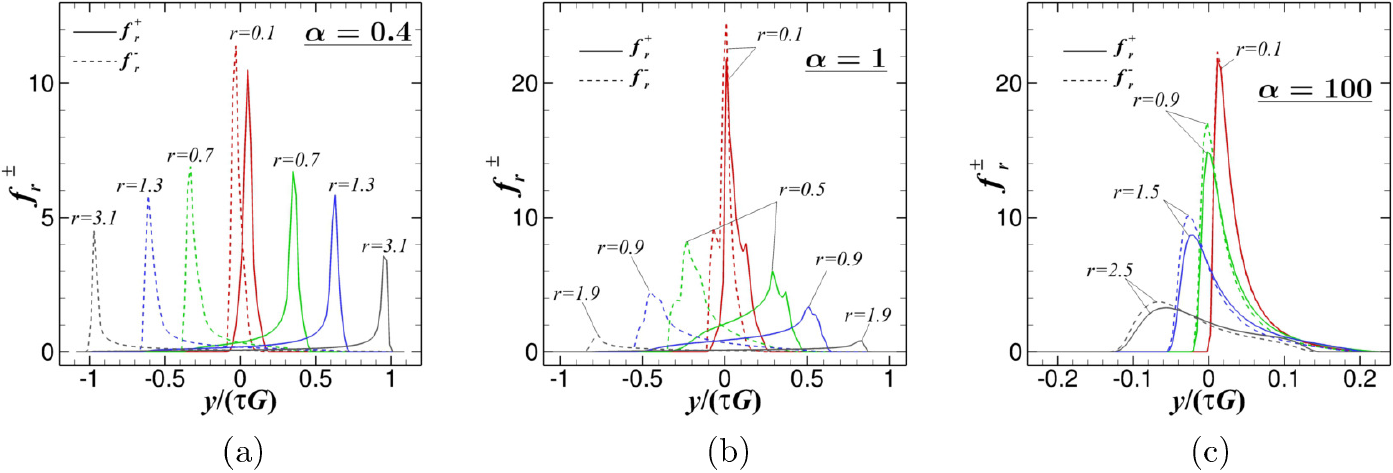}
\caption{
The stationary distributions of the internal state $y$ at different distances from the center of the aggregate $r=|x-x_0|$ at the moderate stiffness $\delta=0.1$.
Figures (a), (b), and (c) show the results at different relative adaptation times, i.e., $\alpha=0.4$, $\alpha=1$, and $\alpha=100$, respectively, at the mean tumbling frequency $\lambda_0=10$.
Here, $G$ denotes the maximum value of the spatial gradient of $M(S)$, i.e., $G=\max_x |\partial_x M(S)|$.
}\label{fig.disty_eps01dlt01}
\end{figure}
\begin{figure}[htbp]
\centering
\includegraphics[width=0.85\textwidth]{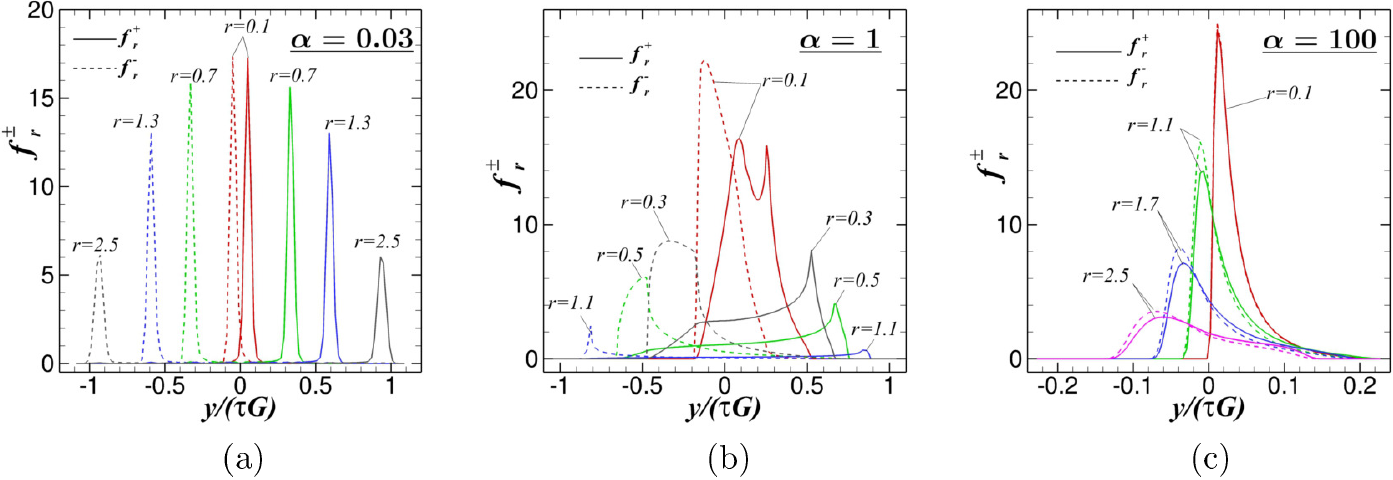}
\caption{
The stationary distributions of the internal state $y$ at different distances from the center of the aggregate $r=|x-x_0|$ at the large stiffness $\delta=0.01$.
Figures (a), (b), and (c) show the results at different relative adaptation times, i.e., $\alpha=0.03$, $\alpha=1$, and $\alpha=100$, respectively, at the mean tumbling frequency $\lambda_0=10$. See also the caption in Fig. \ref{fig.disty_eps01dlt01}.
}\label{fig.disty_eps01dlt001}
\end{figure}

Figures \ref{fig.disty_eps01dlt01} and \ref{fig.disty_eps01dlt001} show the stationary distributions of the internal state $y$ at different distances from the center of the aggregate $r=|x-x_0|$, which are defined as
$$
f_r^\pm(y)=\frac{f^+(x_0\mp r,y)+f^-(x_0\pm r,y)}{2}.
$$
That is, $f_r^\pm(y)$ denotes the local distributions of internal state $y$ for the bacterial moving toward and away from the maximum aggregation density at the distance $r$ from the center of the aggregate.
Figure \ref{fig.disty_eps01dlt01} shows the result at the moderate stiffness $\delta=0.1$, while Fig.~\ref{fig.disty_eps01dlt001} shows the results at the large stiffness $\delta=0.01$.
The parameter sets used in Figs.~\ref{fig.disty_eps01dlt01} and \ref{fig.disty_eps01dlt001} are the same as those in Figs.~\ref{fig.macro.delta01} and \ref{fig.macro.delta001}, respectively.

The distribution of the internal state is highly affected by the relative adaptation time.
When the relative adaptation time is short (i.e., Figs. \ref{fig.disty_eps01dlt01}(a) and \ref{fig.disty_eps01dlt001}(a)), $f_r^+$ and $f_r^-$ are symmetric to each other with respect to $y=0$ and have steep peaks at different values of the internal state according to the distance $r$.
Table \ref{t.peaky} shows the relation between the peak position of $f^\pm_r(y)$ with respect to $y$, $y_p^\pm$, and the local spatial gradient of $M(S)$ at the distance $r$.
The following relation almost holds at each distance $r$, i.e.,
$$
y_p^\pm=\pm\tau |\partial_x M(S)|.
$$

\begin{table}[htbp]
\centering
\begin{tabular}{c ccc c ccc}
\hline\hline
&\multicolumn{3}{c}{Fig. \ref{fig.disty_eps01dlt01}(a)}&&\multicolumn{3}{c}{Fig. \ref{fig.disty_eps01dlt001}(a)}\\
\cline{2-4} \cline{6-8}
$r$ & $y^+_p/(\tau G)$& $y^-_p/(\tau G)$& $|\partial_x M/G|$ & & $y^+_p/(\tau G)$& $y^-_p/(\tau G)$& $|\partial_x M/G|$ \\
\hline
0.1 & 0.05 &-0.03 &0.07 &&0.05 &0.05 &0.05\\
0.7 & 0.35 &-0.33 &0.37 &&0.33 &0.33 &0.35\\
1.3 & 0.63 &-0.61 &0.63 &&0.59 &0.59 &0.62\\
2.5 & 0.97 &-0.95 &0.99 &&0.93 &0.93 &0.97\\
\hline\hline
\end{tabular}
\caption{
Relation between the peak positions of $f^\pm_r(y)$ with respect to $y$, $y_p^\pm$ in Fig. \ref{fig.disty_eps01dlt01}(a) and Fig. \ref{fig.disty_eps01dlt001}(a), and the local spatial gradient of $M(S)$.
See also the caption in Fig. \ref{fig.disty_eps01dlt01}.
}\label{t.peaky}
\end{table}

This behavior is intuitively explained as follows.
From Eq.~(\ref{eq_m}) and the definition of $y$, the dynamics of the internal state $y$ of each bacterium is described as
$$
\dot y=\dot M(S) -\frac{y}{\tau},
$$
where $\dot M(S)$ is the temporal derivative of $M(S)$ along the moving path of each bacterium and is replaced with $|\partial_x M(S)|$ (or $-|\partial_x M(S)|$) at each local position $r$ when the bacteria move toward (or away from) the maximum aggregation density.
Note that $\dot M(S)$ changes the sign due to the tumbling of bacteria.
Thus, when we denote the internal states of the bacteria moving toward (or away from) the maximum aggregation density $y^+$ (or $y^-$), and the temporal evolution of $y^\pm$ at each instant is written as
\begin{equation}\label{eq.internaly}
\dot y^\pm=\pm|\partial_x M| -\frac{y}{\tau}.
\end{equation}
When the adaptation time $\tau$ is much smaller than the run time $\lambda_0^{-1}$, i.e., $\lambda_0\tau\ll1$, the internal state is determined by the local equilibrium state to be $y_p^\pm=\pm\tau|\partial_x M|$ in each run duration.

On the other hand, at large relative adaptation time $\alpha=100$ (i.e., Figs. \ref{fig.disty_eps01dlt01}(a) and \ref{fig.disty_eps01dlt001}(a)), both distributions $f_r^\pm$ concentrate around $y=0$ in the scaled internal variable $y/(\tau G)$, while at the moderate relative adaptation time $\alpha=1$ (i.e., Figs. \ref{fig.disty_eps01dlt01}(b) and \ref{fig.disty_eps01dlt001}(b)), the internal state moderately concentrates around $y_p^\pm=\pm\tau|\partial_x M|$ at each local position $r$.

These behaviors of the internal state variable $y$, according to the change in the relative adaptation time $\alpha$, are consistent with the continuum-limit solutions obtained by the asymptotic analysis, which is presented in the next section and in Ref.~\cite{PSTY2020}; that is,
the continuum limit solution at $\lambda_0^{-1}\rightarrow 0$ is obtained as
$f_0^\pm(t,x,y)=\rho(t,x)\delta(\frac{y}{\tau}=\pm\partial_x M(S))$ when $\alpha\ll 1$ and $f_0^\pm(t,x,y)=\rho(t,x)\delta(y=0)$ when $\alpha\gg 1$.
When $\alpha=1$, the distribution of the internal state moderately concentrates around $y=\pm\tau\partial_x M(S)$.
This characteristic behavior of the internal state variable with respect to the relative adaptation time $\alpha$ is less affected by the change in the stiffness parameter $\delta$ compared to the dependency of the macroscopic population density $\rho$ on the stiffness parameter $\delta$.
In the next subsection, we consider how the stiffness parameter $\delta$ affects the individual motions of bacteria to create different aggregation profiles, as shown in Figs.~\ref{fig.macro.delta01} and \ref{fig.macro.delta001}.

\subsubsection{Local mean run length}

\begin{figure}[htbp]
\centering
\includegraphics[width=0.85\textwidth]{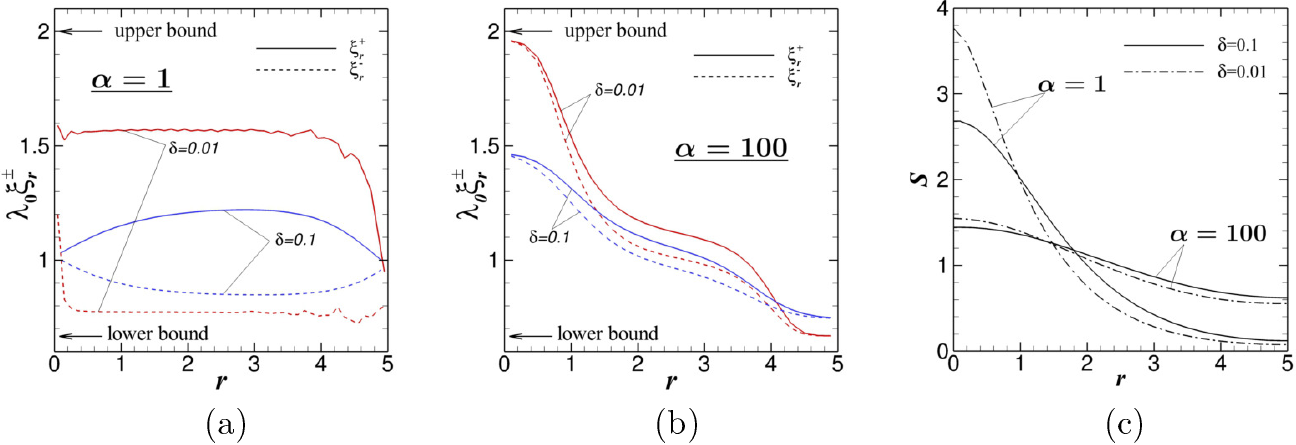}
\caption{
Spatial distributions of the mean run length $\xi^\pm_r$, defined by Eq. (\ref{eq_xi}), and the chemical cue $S$ in the stationary state at different values of the stiffness parameter $\delta$.
Figures (a) and (b) show the results of the mean run length at $\alpha=1$ and $\alpha=100$, respectively, and Figure (c) shows the results of chemical cues at the same parameter sets.
The mean tumbling frequency $\lambda_0=10$ is fixed.
}\label{fig.xi}
\end{figure}
Figure \ref{fig.xi} shows the spatial distributions of the local mean run length of bacteria at moderate stiffness $\delta=0.1$ and at large stiffness $\delta=0.01$.
Here, the local mean run length $\xi_r^\pm$ is calculated as
\begin{equation}\label{eq_xi}
\xi_r^\pm=
\int \frac{f_r^\pm(y)}{\rho_r \lambda_0\Lambda(\frac{y}{\delta})}dy,
\end{equation}
where $\rho_r$ is the population density at distance $r$.
Thus, $\xi_r^+$ denotes the local mean run length of the bacteria moving toward the maximum aggregation density at distance $r$, while $\xi_r^-$ denotes that of the bacteria moving away from the maximum aggregation density at distance $r$.
Since the modulation function $\Lambda_\delta(y)$, which is defined by Eq. (\ref{eq_Lambda}), is bounded, the mean run length is also bounded.
In Figs. \ref{fig.xi} (a) and (b), the upper and lower bounds of the mean run length are shown by leftward arrows on the vertical axis.

At moderate relative adaptation time $\alpha=1$ (Fig.~\ref{fig.xi}(a)), there is a significant difference between $\xi_r^+$ and $\xi_r^-$.
This indicates that the bacteria create highly biased motions according to the moving directions.
At moderate stiffness $\delta=0.1$, the biased motion is maximized around the location of the maximum gradient of $M(S)$, i.e., $|\partial_x M(S)|=|\partial_x S/S|$ (see Fig.~\ref{fig.xi}(c)), while at large stiffness $\delta=0.01$, the biased motion is further enhanced over the whole domain except the vicinities at $r=0$ and $r=5$ so that the spike-like aggregation profile forms due to the highly biased motions of bacteria (see also Fig.~\ref{fig.macro.delta001}(b)).

On the other hand, at large relative adaptation time $\alpha=100$, both $\xi_r^+$ and $\xi_r^-$ vary similarly according to the local amplitude (not the spatial gradient) of the chemical cue $S$, so the biased motion of bacteria is less prominent.
This observation is consistent with the $y$ distributions in Figs.~\ref{fig.disty_eps01dlt01}(c) and \ref{fig.disty_eps01dlt001}(c), where the difference between $f_r^+$ and $f_r^-$ is small.
Thus, aggregation is weakened compared to that at moderate relative adaptation time $\alpha=1$.

In contrast to the results at $\alpha$=1 [Fig.~\ref{fig.xi}(a)], at $\alpha$=100 [Fig.~\ref{fig.xi}(b)], the stiffness of the chemotactic response $\delta$ does not enhance the biased motions in different moving directions but amplifies the spatial modulation of the mean run lengths in both moving directions, although the spatial profile of the chemical cue is less affected by the stiffness parameter [Fig.~\ref{fig.xi}(c)].
Remarkably, at large stiffness $\delta=0.01$, the mean run lengths in both moving directions $\xi^\pm$ attain the upper and lower bounds in the vicinities of $r=$0 and 5, respectively.
Thus, the biased motion is rather hindered due to the boundedness of the tumbling frequency modulation $\Lambda_\delta(y)$ in the vicinities of $r=$0 and 5 such that the plateau regimes of the trapezoidal aggregation profile in Fig.~\ref{fig.macro.delta001}(c) are created.

The overall observations in the local mean run length indicate that the biased motion of bacteria is mostly determined by the relative adaptation time $\alpha$, while the stiffness parameter $\delta$ amplifies only the signal of the internal state in the chemotactic response function.
These orthogonal effects of the relative adaptation time $\alpha$ and the stiffness of the chemotactic response $\delta$ produce the variety of aggregation profiles observed in Figs.~\ref{fig.macro.delta01} and \ref{fig.macro.delta001}.

\section{Asymptotic analysis}\label{sec.asymptotic}
\subsection{Keller-Segel limit}
In this section, we formally derive the asymptotic equations at $\eps=\lambda_0^{-1}\rightarrow 0$ under the diffusive scalings of (\ref{eq_St}) and (\ref{eq_f}), i.e.,
\begin{equation}\label{eq_feps}
\eps\partial_t f^\pm_\eps
\pm \partial_xf^\pm_\eps
+\partial_y\left\{
\left(
\eps \partial_t M(S_\eps)\pm\partial_xM(S_\eps)-\frac{y}{\tau}
\right)f^\pm_\eps
\right\}
=\pm\frac{\Lambda_\delta(y)}{2\eps}(f^-_\eps-f^+_\eps).
\end{equation}
\begin{equation}\label{eq_Seps}
\eps\partial_t S_\eps=D_S\partial_{xx} S_\eps -S_\eps +\rho_\eps,
\end{equation}
We consider three different adaptation time scalings, i.e., (I) $\tau=O(\eps)$ (i.e., $\alpha=O(1)$), (II) $\tau=O(\eps^2)$ (i.e., $\alpha=O(\eps)$), and (III) $\tau=O(1)$.
Here, we set the time scaling parameters as $\sigma$=$\sigma_S$=$\eps$.

In Ref.~\cite{PSTY2020}, asymptotic analysis of the continuous velocity version of Eq.~(\ref{eq_f}) is carried out at the same scalings of the adaptation time by assuming that the stiffness parameter $\delta$ is the same order as the adaptation time, i.e., $\delta=O(\tau)$, and the spatial gradient of $M(S)$ is uniform, i.e., $\nabla_x M(S)=G$, where $G$ is constant.
The results of Ref.~\cite{PSTY2020} are briefly summarized as follows:
\begin{itemize}
\item In Case I, supposing $\chi=O(\eps)$, a hyperbolic model is found at $\sigma=1$, while a novel type of flux-limited KS model is found at $\sigma=\eps$. Furthermore, in both time scalings ($\sigma$=1 and $\eps$), the leading order solution $f_0$ is obtained from a new type of equilibrium equation (which is the continuous velocity version of Eq.~(\ref{apeq_Q0})).
\item In Case II, supposing $\chi=O(\eps)$, a FLKS model is obtained at $\sigma=\eps$, and the leading order solution is explicitly written as $f_0=\rho(t,x)\delta(y=v\cdot G)$.
\item In Case III, a classical KS-type model is obtained for $\sigma=\eps$, and the leading order solution is explicitly written as $f_0=\rho(t,x)\delta(y=0)$.
\end{itemize}

In this paper, we consider the case where both stiffness and modulation are moderate, i.e., $\delta=O(1)$ and $\chi=O(1)$.
By taking the sum of Eq. (\ref{eq_feps}), we have
$$
\partial_t\left(\frac{f_\eps^++f_\eps^-}{2}\right)
+\partial_x\left(\frac{f_\eps^+-f_\eps^-}{2\eps}\right)
+\partial_y\left\{
\left(\eps\partial_tM_\eps-\frac{y}{\tau}\right)
\left(\frac{f_\eps^++f_\eps^-}{2}\right)
+\partial_x M_\eps \left(\frac{f_\eps^+-f_\eps^-}{2\eps}\right)
\right\}=0.
$$
Integration of the above equation with respect to $y$ gives the following macroscopic conservation law:
\begin{equation}\label{eq_rhoeps}
\partial_t\rho_\eps+\partial_x \left(\frac{J_\eps}{\eps}\right)=0,
\end{equation}
where the flux $J_\eps$ is defined as
\begin{equation}\label{eq_jeps}
J_\eps=\int_R \frac{f_\eps^+-f_\eps^-}{2}dy.
\end{equation}
The continuum-limit equations for the population density $\rho_0$ can be derived from the above formulas (\ref{eq_rhoeps}) and (\ref{eq_jeps}).

By taking the continuum limit at Eq. (\ref{eq_Seps}), we also obtain that $S_0$ is the following equation for chemoattractant $S_0$:
\begin{equation}\label{eq_S0}
-D_S\partial_{xx}S_0+S_0=\rho_0.
\end{equation}

We carried out asymptotic analysis under different scalings of the adaptation time, i.e., (i) $\tau=O(\eps)$, (ii) $\tau=O(\eps^2)$, and (iii) $\tau=O(1)$, and obtained the KS-type models as described below.
The detailed calculations are given in \ref{app_asymp_1}--\ref{app_asymp_3}.
Here, we only summarize the main results, i.e., for (i) $\tau=O(\eps)$ (or $\alpha=O(1)$),
\begin{equation}\label{eq_c1_KS}
\partial_t\rho_0-\partial_{xx}\rho_0+\partial_x\left(\frac{\alpha \chi \partial_x M(S_0)}{\delta(1+\alpha)}\rho_0\right)=0,
\end{equation}
for (ii) $\tau=O(\eps^2)$ (or $\alpha=O(\eps)$),
\begin{equation}\label{eq_c2_KS}
\partial_t \rho_0-\partial_{xx}\rho_0=0,
\end{equation}
and for (iii) $\tau=O(1)$ (or $\alpha=O(1/\eps)$),
\begin{equation}\label{eq_c3_KS}
\partial_t\rho_0-\partial_{xx}\rho_0+\partial_x\left(\frac{\chi\partial_x M_0}{\delta}\rho_0\right)=0.
\end{equation}
Note that Eqs.~(\ref{eq_c2_KS}) and (\ref{eq_c3_KS}) coincide with Eq.~(\ref{eq_c1_KS}) at the limits $\alpha\rightarrow 0$ and $\alpha \rightarrow \infty$, respectively.
Thus, the KS-type equation (\ref{eq_c1_KS}) can uniformly describe the continuum-limit behavior of the kinetic transport equation when the adaptation time is at most moderate, i.e., $\tau < O(1)$.

\subsubsection{Linear instability of the KS system}\label{sec.ksinstability}
The linear instability of the KS system around the uniform solution $\rho=S=1$ is obtained as follows:

First, we consider a small perturbation of the uniform solution in the following form:
$$
\rho(t,x)=1+\tilde \rho(x)e^{\mu t},\quad S(t,x)=1+\tilde S(x)e^{\mu t},
$$
and linearize Eq.~(\ref{eq_c1_KS}) as
\begin{eqnarray*}
&\mu \tilde \rho(x)e^{\mu t}-\tilde \rho''(x)e^{\mu t}
+\partial_x\left(
\frac{\alpha \Lambda'_\delta(0)}{1+\alpha}\tilde S'(x) e^{\mu t}(1+\tilde \rho'(x) e^{\mu t}
\right)=0,\\
&\mu \tilde \rho(x)-\tilde \rho''(x)
+\left(
\frac{\alpha \Lambda'_\delta(0)}{1+\alpha}\right)\tilde S''(x)=0.
\end{eqnarray*}
By taking the Fourier transform of the above equations, we obtain
\begin{equation}\label{eq_rhok}
\mu \tilde \rho_k+k^2\left(
\tilde \rho_k-
\frac{\alpha \Lambda'_\delta(0)}{1+\alpha} \tilde S_k
\right)=0,
\end{equation}
where $k$ is the Fourier variable and $\tilde \rho_k$ and $\tilde S_k$ are the Fourier transforms of $\tilde \rho(x)$ and $\tilde S(x)$, respectively, which are calculated as $\tilde \rho_k=\int_R \tilde \rho(x)e^{-\mathrm{i} k x}dx$.
By inserting the Fourier transform of Eq.~(\ref{eq_S0}),
\begin{eqnarray*}
(1+D_Sk^2)\tilde S_k=\tilde \rho_k,\\
\tilde S_k=\frac{\tilde \rho_k}{1+D_Sk^2},
\end{eqnarray*}
into Eq. (\ref{eq_rhok}), we obtain
\begin{equation}
\left[
\mu+k^2\left(
1-\frac{\alpha \Lambda'_\delta(0)}{(1+\alpha)(1+D_Sk^2)}
\right)
\right]\tilde \rho_k=0	.
\end{equation}

Hence, the condition in which the mode $k$ becomes linearly unstable ($\mu>0$) is written as
\begin{equation}\label{eq_ksinstability}
\Lambda'_\delta(0)>\frac{1+\alpha}{\alpha}(1+D_Sk^2).
\end{equation}

Thus, the instability of the mode $k$ occurs when the stiffness of the chemotactic response $\Lambda'_\delta(0)$ is larger than the right-hand side of Eq.~(\ref{eq_ksinstability}).
This also indicates that when the stiffness of the response $\Lambda'_\delta(0)$ and the diffusion coefficient $D_S$ are fixed, instability more likely occurs as $\alpha$ increases.

\subsection{A novel asymptotic equation at large adaptation time}
In the previous section, we show that the KS-type model involving the relative adaptation time $\alpha$, Eq.~(\ref{eq_c1_KS}) is derived when the adaptation time is scaled as $\tau=O(\eps^n)$ ($n=0,1,2$) at the continuum limit $\eps=\lambda_0^{-1}\rightarrow 0$.

In this section, we consider the case where the adaptation time is very large, i.e., $\tau=O(\eps^{-1})$.
The formal asymptotic analysis at the large adaptation-time regime, i.e., $\tau=\tilde \tau/\eps$ with $\tilde \tau=O(1)$, is carried out in \ref{app_asymp_largetau}, and the following asymptotic equation is obtained at the limit $\eps\rightarrow 0$:
\begin{equation}\label{eq_limitp}
    \partial_t p_0
    -\partial_x\left(\frac{\partial_x p_0}{\Lambda_\delta(M(S)-m)}\right)
    +\partial_m\left(
    \frac{M(S)-m}{\tilde \tau}p_0
    \right)=0,
\end{equation}
where $p_0=p_0(t,x,m)$ is the continuum-limit solution to Eq.~(\ref{eq_p}) at the large adaptation-time scaling $\tau=\tilde \tau/\eps$.
This novel asymptotic equation retains the internal state variable $m$ as an independent variable. 
Hence, the population density $\rho_0$ is obtained by the integration of $p_0$ with respect to the internal variable $m$:
\begin{equation}
\rho_0(t,x)=\int_{-\infty}^\infty p_0(t,x,m)dm.
\end{equation}

Obviously, the novel asymptotic equation is completely different from the KS system.
However, we can confirm the consistency of this asymptotic equation with the KS system (\ref{eq_c3_KS}) at the limit $\tilde{\tau}\rightarrow 0$ (see also \ref{app_asymp_largetau}).
We will also numerically confirm the robustness of the asymptotic solution to Eq.~(\ref{eq_limitp}) in the next section.

Since $\tau$ and $\lambda_0$ are nondimensionalized as Eq.~(\ref{def_parap}), the large adaptation-time scaling $\tau\sim \lambda_0$ is rewritten in dimensional form as
\begin{equation}\label{eq_td}
\tau\sim t_d, \quad t_d=\frac{L_0^2}{D_\rho},
\end{equation}
where $t_d$ is the time for bacterial population to diffuse over the characteristic length $L_0$, which is, in the present problem, determined by the diffusion of chemoattractant in the medium (Eq.~(\ref{eq_L0}), and $D_\rho$ is the diffusion constant of the macroscopic population density defined as $D_\rho=V_0^2/\lambda_0$.
Thus, the novel asymptotic equation (\ref{eq_limitp}) is appropriate when the adaptation time is comparable to the diffusion time of the population density in the characteristic length, $\tau=O(t_d)$, while the KS system is only valid when the adaptation time is much smaller than the diffusion time, $\tau\ll t_d$.

\subsection{Asymptotic behavior at a large adaptation time}\label{sec.discussion}
In this section, we further discuss the asymptotic behavior of the population density with respect to $\eps=\lambda_0^{-1}$ at large adaptation times, i.e., (i) $\tau=O(1)$ and (ii) $\tau=O(\eps^{-1})$.
\begin{figure}[htbp]
\centering
\includegraphics[width=0.85\textwidth]{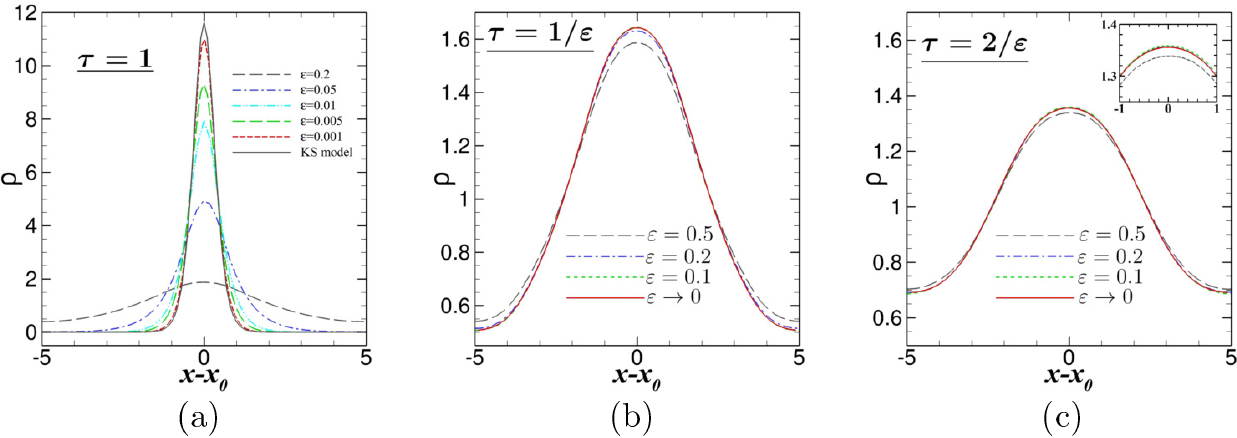}
\caption{
Asymptotic behaviors of the population density $\rho$ at moderate stiffness $\delta=0.1$.
Figures (a), (b), and (c) show the results at different scalings of the adaptation time, i.e., $\tau=1$, $\tau=1/\eps$, and $\tau=2/\eps$, respectively, with $\eps=\lambda_0^{-1}$.
The inset in (c) shows the magnification around the center of the aggregate.
The solid line in (a) shows the result of the KS model, while the solid lines in (b) and (c) show the results of the novel asymptotic equation (\ref{asymp_largetau}).
}\label{fig.asymptotic_del01}
\end{figure}
Figure \ref{fig.asymptotic_del01} shows the asymptotic behaviors of the population density at moderate stiffness $\delta=0.1$.
At $\tau=1$, the MC results approach those of the KS system as $\eps$ decreases.
This observation is consistent with the asymptotic analysis in Sec.~\ref{sec.asymptotic}.
However, the asymptotic convergence is very slow; a significant deviation remains between MC results and the KS system results, even at small $\eps$, e.g., $\eps\sim 0.01$.

This slow asymptotic convergence is also confirmed in Fig.~\ref{fig.maxrho}, where the MC results at $\tau=1$ line up in an upper-right direction and gradually approach the KS limit as the relative adaptation time $\alpha$ increases.
On the other hand, it is also seen that when $\alpha$ is fixed at $\alpha=O(1)$, i.e., $\tau=O(\eps)$, the MC results converge to the KS result more rapidly as $\eps$ decreases.

\begin{figure}[htbp]
\centering
\includegraphics[width=0.85\textwidth]{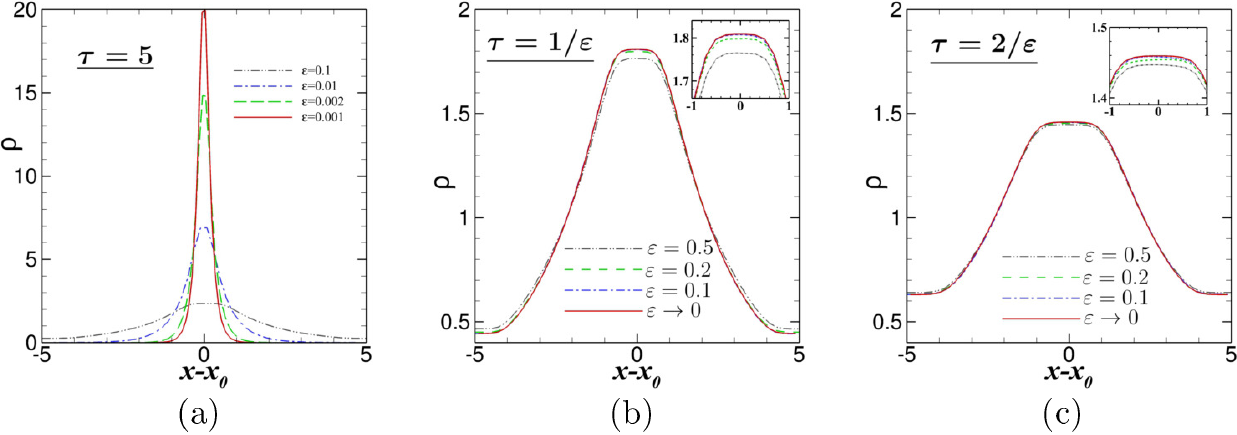}
\caption{
Asymptotic behaviors of the population density $\rho$ at large stiffness $\delta=0.01$.
Figures (a), (b), and (c) show the results at different scalings of the adaptation time, i.e., $\tau=5$, $\tau=1/\eps$, and $\tau=2/\eps$, respectively.
The solid lines in (b) and (c) show the results of the novel asymptotic equation (\ref{asymp_largetau}).
}\label{fig.asymptotic_del001}
\end{figure}

When the adaptation time is set to $\tau=\tilde \tau/\eps$ (e.g., $\tilde \tau$=1 in Fig.~\ref{fig.asymptotic_del01}(b) and $\tilde \tau$=2 in Fig.~\ref{fig.asymptotic_del01}(c)), the asymptotic convergence of the MC results is much faster at $\tau=O(\eps^{-1})$ compared to that observed at $\tau=1$.
Remarkably, even at moderately small $\eps$, e.g., $\eps\sim 0.1$, the MC results almost coincide with the numerical solutions of the novel asymptotic equation (\ref{eq_limitp}).
Convergence to the asymptotic solution is also observed in Fig.~\ref{fig.asymptotic_del001}.
These results numerically confirm that the asymptotic solution to Eq.~(\ref{eq_limitp}) is robust at the large adaptation-time scaling $\tau=O(\eps^{-1})$.

Note that the maximum aggregation density decreases as $\tilde \tau$ increases, as is already observed in Fig.~\ref{fig.maxrho}.
Thus, the nonmonotonic behavior of the maximum aggregation density with respect to the adaptation time $\tau$ at each fixed $\eps$ can be viewed as the transition from the KS-type solution to the asymptotic solution to Eq.~(\ref{eq_limitp}).

Figure \ref{fig.asymptotic_del001} shows the asymptotic behaviors of the population density at large stiffness $\delta=0.01$.
It is seen from Fig. \ref{fig.asymptotic_del001}(a) that when $\tau=5$ is fixed, the width of the aggregate narrows as $\eps$ decreases, so  asymptotic convergence is not confirmed from the present MC results.

Interestingly, a trapezoidal aggregate is robustly formed at the large adaptation-time scaling $\tau=O(1/\eps)$; in Figs.~\ref{fig.asymptotic_del001}(b) and (c), the MC results at moderately small $\eps$, e.g., $\eps\lesssim 0.1$, are close to each other and asymptotically converge to the trapezoidal profile.
Furthermore, it is clear  that the trapezoidal profile is related to the novel asymptotic equation (\ref{eq_limitp}) at the large adaptation-time scaling $\tau=\tilde \tau/\eps$.

\section{Summary and perspectives}\label{sec.summary}
We investigated the self-organized aggregation of chemotactic bacteria in one-dimensional space with periodic boundary conditions based on a two-stream kinetic transport model with an internal state coupled with the chemoattractant equation.
MC simulations were conducted for a wide range of adaptation times $\tau$ at various values of the mean tumbling frequency $\lambda_0$ and the stiffness of chemotactic response $\delta$
Asymptotic analysis of the kinetic transport model was also carried out to complement the MC results.
Thus, the effect of the adaptation time on the aggregation behavior was investigated both macroscopically and microscopically.

An important finding is the nonmonotonic dependence of the adaptation time on the aggregation behavior.
See, for example, Figs.~\ref{fig.instability} and \ref{fig.maxrho}.
A sharp transition between stable and unstable modes is observed when the relative adaptation time $\alpha$, which is defined as $\alpha=\lambda_0\tau$, is small, e.g., $\alpha \lesssim 1$.
In this small $\alpha$ regime, instability always occurs when $\alpha$ is slightly larger than the critical value of the linear instability condition of the KS model, which is derived from the kinetic transport model by asymptotic analysis, and the maximum aggregation density rapidly increases as $\alpha $ increases.
However, when the relative adaptation time is large, e.g., $\alpha \gtrsim \lambda_0$, the MC results deviate from the linear stability condition of the KS model, and the maximum aggregation density gradually decreases as $\alpha$ increases.
Thus, there exists an optimal adaptation time to enhance aggregation around $1 < \alpha < \lambda_0$.
This nonmonotonic behavior is a significantly important feature that can be described at the kinetic level but not at the KS level; in the KS model, the maximum aggregation density monotonically increases and saturates to a certain value as $\alpha$ increases.

We also investigated the microscopic behaviors in the variety of aggregation profiles in terms of the local distribution of the internal state (Figs.~\ref{fig.disty_eps01dlt01} and \ref{fig.disty_eps01dlt001}) and the spatial distribution of the local mean run length (Fig.~\ref{fig.xi}).
We found an orthogonal effect of the relative adaptation time $\alpha$ and the stiffness of the chemotactic response $\delta$ on the microscopic dynamics.
That is, the relative adaptation time $\alpha$ significantly affects the distribution of the internal state, while the stiffness $\delta$ does not similarly affect the distribution of the internal state but only amplifies the signal of the internal state in the response function $R_\delta(y)$.

From these microscopic perspectives, the optimal adaptation time to produce sharp aggregation can be intuitively explained as follows.
When the adaptation time $\tau$ is smaller than the mean run duration $\lambda_0^{-1}$ (i.e., $\alpha<1$), the internal state $y$ is rapidly equilibrated at $y_p^\pm=\pm\tau|\partial_x M|$ depending on the moving direction in each run duration.
Here, since the equilibrium state $M$ sensed by the bacteria is temporally changed along the run of each bacterium, the internal state $y$ is not equilibrated at $y=0$ but is equilibrated at $y=y_p$, which is linearly proportional to the adaptation time $\tau$ unless the bacteria tumbles.
Thus, the amplitude of the chemotactic response $|R(y)|$ (see Eq.~(\ref{eq_Lambda})) becomes larger in each run duration at adaptation time $\tau$ when $\alpha<1$.
On the other hand, when the adaptation time $\tau$ is much larger than the run duration $\lambda_0^{-1}$ (i.e., $\alpha >\lambda_0 \gg 1$), the internal state cannot be equilibrated in each run duration so that the biased motion between the different moving directions is significantly reduced, as shown in Fig.~\ref{fig.xi}(b).
Thus, aggregation is hindered when $\alpha$ is very large.
These competitive effects of the adaptation time on the chemotactic response lead to the optimal behavior of chemotactic aggregation.

\begin{figure}
\centering
\includegraphics[width=0.4\textwidth]{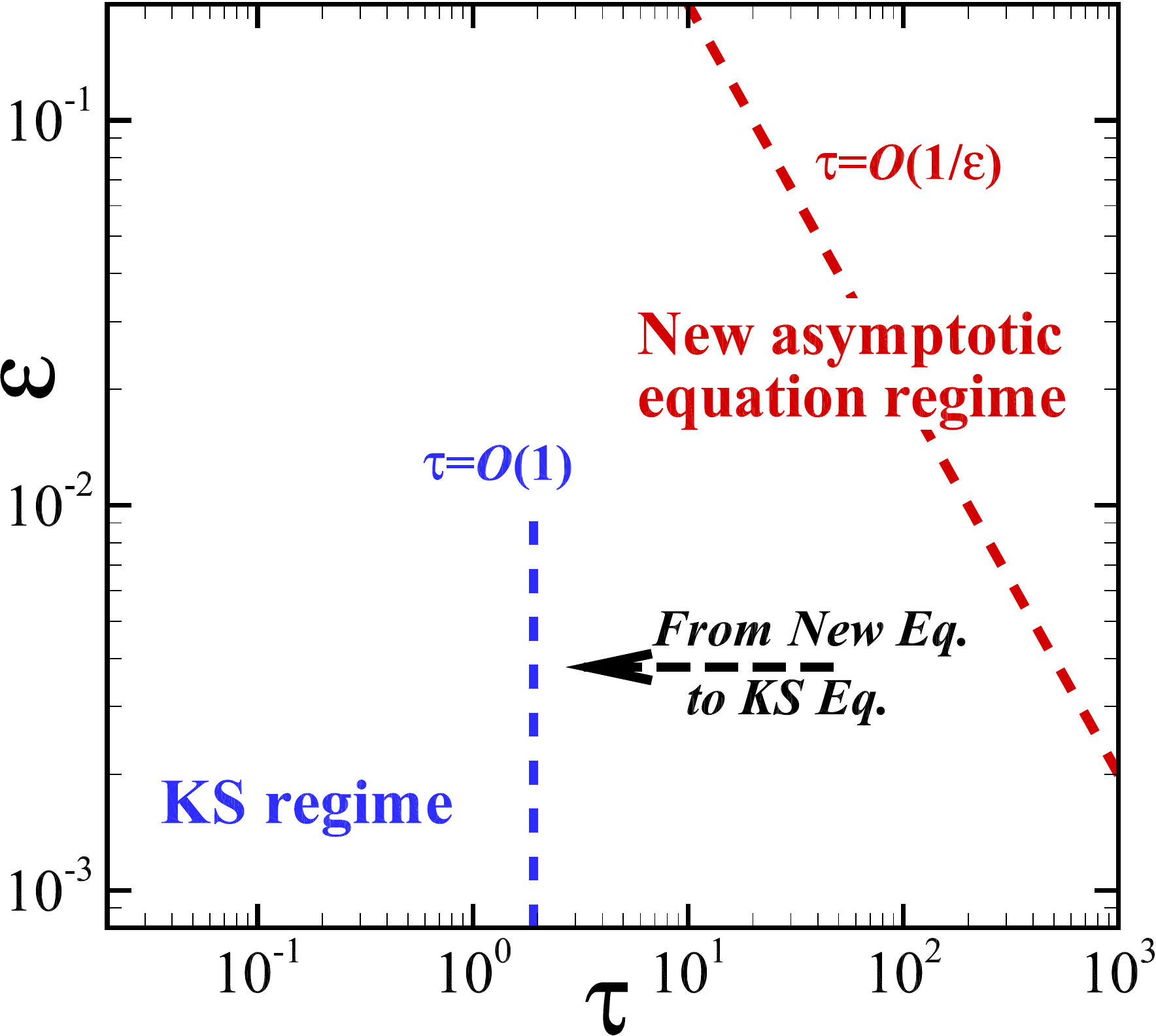}
\caption{Schematic of the asymptotic regimes.
The novel asymptotic equation (\ref{eq_limitp}) is discovered along the red solid line while the KS regime is limited in $\tau\lesssim O(1)$.
The KS equation is also obtained from the novel asymptotic equation by taking a consistent limit.
}\label{fig.cartoon}
\end{figure}
Another remarkable result of this paper is the discovery of the novel asymptotic equation Eq.~(\ref{eq_limitp}) at large adaptation-time scaling $\tau=O(\eps^{-1})$.
This asymptotic regime physically indicates that the adaptation time $\tau$ is comparable to the diffusion time of the population density in the characteristic length [Eq.~(\ref{eq_td})].
A numerical comparison clarified that the trapezoidal aggregates, which are robustly formed at the large adaptation-time regime are well described by the novel asymptotic equation. See Fig.~\ref{fig.asymptotic_del001}.

Figure \ref{fig.cartoon} is the cartoon of the asymptotic regimes of the novel asymptotic equation (\ref{eq_limitp}) and the KS system (\ref{eq_c1_KS}).
The KS system is not valid even near the continuum limit $\eps\ll 1$ unless the adaptation time $\tau$ is at most moderate $\tau<O(1)$.
Although the novel asymptotic equation is obtained at the scaling $\tau=O(\eps^{-1})$ (the red dashed line in the figure), we also formally show that the novel asymptotic equation converges to the KS system when $\tilde{\tau}(=\eps\tau)\rightarrow 0$ (the dashed left arrow in the figure).

The asymptotic behavior of the MC results clearly illustrates the limitation of the appropriate adaptation-time regime for the KS system and the transient behavior from the KS solution to the novel asymptotic solution when the adaptation time becomes large.
Thus, the nonmonotonic behavior of the aggregation density with respect to the adaptation time can be interpreted as the transient behavior from the KS-type solution to the novel asymptotic solution, although the detailed analysis of the transient behavior remains as important future work.

One may think that such a large adaptation time is biologically unrealistic.
However, when the system size is as small as $L_0\sim 100\, \mu\mathrm{m}$ and the mean run length of bacteria is measured as $l_0\sim 20\,\mu\mathrm{m}$, the ratio of the mean run length of bacteria to the system size, $\eps$, is typically estimated to be $\eps\sim 0.2$.
In this case, the relative adaptation time in the novel asymptotic regime is estimated to be $\alpha\sim 25$.
The adaptation time of bacteria is usually much longer than the run time, so $\alpha\sim 25$ is not unrealistic but rather commonly observed.

Indeed, when the system size is as small as $\eps \sim 0.1$, unusual aggregation behaviors of bacteria are observed in experiments. For example, in Ref.~\cite{MBBO2003}, the volcano-like aggregation profile of $\textit{E. coli}$ is observed at $L_0\sim 100\,\mu\mathrm{m}$, and in Ref.~\cite{BM2003}, swarm bands of marine bacteria around a chemoattractant microbead with a distance of approximately $L_0\sim 20\,\mu\mathrm{m}$ are observed.
The common feature in these curious aggregation behaviors is that the maximum aggregation density is not located at the center of the aggregate but rather at a certain distance from the center.
It is interesting that the trapezoidal aggregation profile is obtained in the present MC simulation in the same parameter regime as in the experiments, although neither volcano-like aggregation nor swarm rings have yet been confirmed.

In the present paper, to focus on the multiscale mechanism between collective motion and individual motility involving internal adaptation dynamics, we ignore several factors that should be important in reproducing experimental results.
Among them, the nutrients consumed by bacteria and proliferation due to cell division should play a significant role in the collective motions and pattern formations of chemotactic bacteria.
For example, in Ref.~\cite{CMD2013}, the variety of patterns observed in experiments was successively reproduced by an individual-based simulation coupled with the consumption of nutrients and the secretion of chemoattractants.
Inclusion of these factors in the present kinetic transport model is rather straightforward and plays a significant role in collective motions and pattern formations.
Indeed, in our previous study~\cite{PY2018}, the instability that leads to Turing-like periodic pattern formation was clarified based on the kinetic transport model involving the proliferation of bacteria.
The traveling pulse created by the bacteria pursuing the nutrient is also successively illustrated based on the kinetic transport model coupled with the reaction-diffusion equations of the nutrient and secreted chemoattractant~\cites{SCBPBS2011,C2020}.
In these previous kinetic studies, the internal dynamics were completely ignored.
Typically, the adaptation time is comparable to the characteristic time of the collective traveling pulse of bacteria pursuing the nutrient; for example, in the traveling pulse observed in Ref.~\cite{SCBPBS2011}, the characteristic time of the wave is estimated to be $t_0\sim 20$ s.
Thus, we can expect some coupling effects between the internal adaptation dynamics and the collective traveling pulse pursuing the nutrient.
Investigation of the effect of the adaptation time on the traveling waves and pattern formations coupled with the sensing of the nutrient and secreted chemoattractant is  important future work.

In the present kinetic transport model, we also utilize the simplified model for intracellular signal transduction, where only the adaptation dynamics of the internal state are considered while the excitation dynamics, which are much faster than the adaptation dynamics, are ignored.
The tumbling time (i.e., the time required for the bacteria to change their moving direction) is also ignored since it is much shorter than the run time.
However, in small systems, these fast-time-scale dynamics may bring about a delay in individual motion and affect collective dynamics.
Indeed, in the literature~\cites{BLL2007,JJP2018,SM2011}, the nonunimodal aggregates around a chemoattractant point source are numerically reproduced by using a model involving either the excitation dynamics of the internal state or the tumbling time.

It has also been reported that the noise arising in the intracellular signal transduction process plays a significant role in the occurrence of chemotactic aggregation around a chemoattractant point source~\cite{BCHYLS2019} and in the fractional diffusion mode of chemotactic bacteria~\cites{KEVSC2004,PST2018}.
Furthermore, in Refs.~\cites{MJD2009,MMCD2011}, the origin of the noise in intracellular signal transduction and its effect on the sensing behavior are also argued based on an optimal biochemical network model proposed in Ref.~\cite{KLBTS2005}.
To uncover the multiscale mechanism between collective motion and individual motility involving internal dynamics in a variety of aggregation behaviors, further investigations based on the kinetic transport model involving more sophisticated formulas of the intracellular signal transduction pathway as well as the tumbling state will also be important future work.

\appendix
\section{Derivation of continuum-limit equations}\label{Appendix}

\subsection{Fast adaptation}\label{app_asymp_1}
Here, we consider the case where the adaptation time is comparable to the run duration, i.e., $\tau=O(\eps)$. Thus, the relative adaptation time is $\alpha=O(1)$.
By changing the variable as $g_\eps(t,x,z)=f_\eps(t,x,\tau z)$, Eq.~(\ref{eq_feps}) is written as
\begin{equation}\label{apeq_cI}
\eps^2\partial_t g_\eps^\pm
\pm\eps\partial_x g_\eps^\pm
+\frac{1}{\alpha}\partial_z
\left\{
\left(
\eps\partial_t M_\eps
\pm\partial_x M_\eps-z
\right)g_\eps^\pm
\right\}
=\pm\frac{\Lambda_\delta(\eps \alpha z)}{2}(g_\eps^--g_\eps^+).
\end{equation}
Note that $\Lambda_\delta(\eps\alpha z)=1-R_\delta(\eps\alpha z)$ and $|R_\delta(\eps\alpha z)|\le \frac{\eps \alpha |z|}{\delta}$.
Thus, from the leading term, we obtain
$$
\partial_z[(\pm G_0-z)g_0^\pm]=\pm\frac{\alpha}{2}(g_0^--g_0^+),
$$
where we write $G_0=\partial_x M_0$.
Here, we assume that $g_\eps^\pm$ is compactly supported with respect to $z$.
(This can be proved when $|G_0|$ is bounded, as is done in Sec.~3 of Ref.~\cite{PSTY2020}.)
By integrating the above equation with respect to $z$, we obtain $J_0=0$.

We seek the leading-order solution in the form of $g_0^\pm=\rho_0(t,x) Q_0^\pm(z;G_0)$, where $Q_0^\pm$ is described as
\begin{equation}\label{apeq_Q0}
\partial_z[(\pm G_0-z)Q_0^\pm]=\pm\frac{\alpha}{2}(Q_0^--Q_0^+),
\end{equation}
with $\int_R Q_0^\pm dz=1$.
Furthermore, $Q_0^\pm$ is compactly supported on $z=[-|G_0|,|G_0|]$.

By using the leading order solution, we can write Eq.~(\ref{apeq_cI}) as
\begin{align*}
\pm\partial_x(\rho_0 Q_0^\pm)
&+\frac{1}{\alpha}\partial_z\left\{
\partial_tM_0 \rho_0 Q_0^\pm+(\pm G_0-z)g_1^\pm
\right\}\\
&=\pm \frac{1}{2}(g_1^--g_1^+)\pm \frac{\alpha\chi \rho_0}{2\delta}z (Q_0^+-Q_0^-) + O(\eps),
\end{align*}
where we use $R_\delta(\eps\alpha z)=\frac{\eps\alpha z}{\delta}+O(\eps^2)$ for $z\in[-|G_0|,|G_0|]$.
By integrating the above equation with respect to $z$ and taking the limit $\eps\rightarrow 0$, we obtain
$$
J_1=-\partial_x\rho_0+\frac{\alpha\chi \rho_0}{2\delta}\int_{-|G|}^{|G|}z(Q_0^+-Q_0^-)dz.
$$
Here, we also assume $\int_R |z(g_1^--g_1^+)|dz<+\infty$.
The last term of the above equation is calculated as follows. By integrating Eq. (\ref{apeq_Q0}) multiplied by $z$, we obtain
$$
\int z\partial_z[(\pm G -z)Q_0^\pm]dz
=\pm \frac{\alpha}{2}\int z (Q_0^--Q_0^+),
$$
$$
-\int (\pm G -z)Q_0^\pm dz
=\pm \frac{\alpha}{2}\int z (Q_0^--Q_0^+),
$$
$$
\mp G + \int z Q_0^\pm dz
=\pm \frac{\alpha}{2}\int z (Q_0^--Q_0^+),
$$
$$
\int z (Q_0^+-Q_0^- )dz=\frac{2G}{1+\alpha}.
$$
Thus, we obtain
$$
J_1=-\partial_x \rho_0+\frac{\alpha \chi G \rho_0}{\delta (1+\alpha)}.
$$

Hence, from Eq. (\ref{eq_rhoeps}), we obtain the following KS equation in the continuum limit $\eps\rightarrow 0$,
\begin{equation}\label{apeq_c1_KS}
\partial_t\rho_0-\partial_{xx}\rho_0+\partial_x\left(\frac{\alpha \chi \partial_x M(S_0)}{\delta(1+\alpha)}\rho_0\right)=0.
\end{equation}

\subsection{Very fast adaptation}\label{app_asymp_2}
Here, we consider the case where the adaptation time is much smaller than the run duration, i.e., $\tau=O(\eps^2)$. Hence, $\alpha=O(\eps)$.
By setting $\alpha=\alpha_1\eps$ at Eq. (\ref{apeq_cI}), we obtain
\begin{equation}\label{apeq_cII}
\eps^2\partial_t g_\eps^\pm
\pm\eps\partial_x g_\eps^\pm
+\frac{1}{\alpha_1\eps}\partial_z
\left\{
\left(
\eps\partial_t M_\eps
\pm\partial_x M_\eps-z
\right)g_\eps^\pm
\right\}
=\pm\frac{\Lambda_\delta(\eps^2 \alpha_1 z)}{2}(g_\eps^--g_\eps^+).
\end{equation}
From the leading term of the above equation, we obtain
$$
\partial_z[(\pm\partial_x M_0-z)g_0^\pm]=0.
$$
Thus, the leading order solution is written as
\begin{equation}
g_0^\pm=\rho_0\delta(z=\pm\partial_x M_0),
\end{equation}
where $\delta(z)$ is the Dirac delta function.
The leading order flux is obtained as $J_0=0$.

By taking the difference of Eq. (\ref{apeq_cII}) and integrating it with respect to $z$, we obtain
$$
\eps \partial_x\rho_\eps =-J_\eps +\int_R \frac{R_\delta(\eps^2\alpha_1 z)}{2}(g_\eps^--g_\eps^+)dz+O(\eps^2).
$$
By taking the limit $\eps\rightarrow 0$ of the above equation under the assumption $\int_R|z(g_\eps^--g_\eps^+)|dz < +\infty$, we obtain
$$
J_1=-\partial_x\rho_0.
$$

Hence, we obtain the following diffusion equation at the continuum limit $\eps\rightarrow 0$:
\begin{equation}
\partial_t \rho_0-\partial_{xx}\rho_0=0.
\end{equation}
%

\subsection{Moderate adaptation}\label{app_asymp_3}
Here, we consider the case where the adaptation time is order unity, $\tau=O(1)$.
Hence, $\alpha=O(1/\eps)$.
The leading term of Eq.~(\ref{eq_feps}) gives us
\begin{equation}
f_0^+=f_0^-=f_0,
\end{equation}
and hence, $J_0=0$.

Thus, the next order term is written as
\begin{equation}\label{apeq_c3_first}
\pm\partial_x f_0
+\partial_y\left\{
\left(\pm\partial_x M_0-\frac{y}{\tau}\right)
f_0
\right\}
=\pm\frac{\Lambda_\delta(y)}{2}(f_1^--f_1^+).
\end{equation}
By taking the sum of the above equation, we obtain
$$
\partial_y(y f_0)=0.
$$
Hence, the leading order solution is written as
\begin{equation}\label{apeq_c3_f0}
f_0=\rho_0(t,x)\delta(y=0).
\end{equation}

Thus, from Eq. (\ref{apeq_c3_first}) with Eq. (\ref{apeq_c3_f0}), the flux $J_1$ is calculated as
\begin{eqnarray*}
J_1&=-\int \frac{\delta(y=0)}{\Lambda_\delta(y)}dy\partial_x\rho_0
-\rho_0\partial_x M_0\int \frac{\delta'(y=0)}{\Lambda_\delta(y)}dy,\\
&=-\frac{1}{\Lambda_\delta(0)}\partial_x\rho_0
-\frac{\Lambda_\delta'(0)}{\Lambda_\delta^2(0)}\rho_0\partial_x M_0.
\end{eqnarray*}
Since $\Lambda_\delta(0)=1$ and $\Lambda_\delta'(0)=-\frac{\chi}{\delta}$, we obtain
\begin{equation}\label{apeq_c3_J1}
J_1=-\partial_x\rho_0+\frac{\chi \rho_0\partial_x M_0}{\delta},
\end{equation}
and hence,
\begin{equation}\label{apeq_c3_KS}
\partial_t\rho_0-\partial_{xx}\rho_0+\partial_x\left(\frac{\chi\partial_x M_0}{\delta}\rho_0\right)=0.
\end{equation}

\section{Large adaptation time}\label{app_asymp_largetau}

We consider the kinetic transport equation with internal state $m$, i.e., Eq.~(\ref{eq_p}) under scaling at a very large adaptation time $\tau=\tilde \tau/\eps$,
\begin{equation}\label{ap_peps}
\varepsilon \partial_t p^\pm_\varepsilon
\pm \partial_x p^\pm_\varepsilon
+\varepsilon\partial_m\left(
\frac{M(S)-m}{\tilde \tau}p^\pm_\varepsilon\right)
=\pm\frac{\Lambda_\delta(M-m)}{2\varepsilon}(p_\varepsilon^--p_\varepsilon^+).
\end{equation}

By taking the limit $\varepsilon\rightarrow 0$ in Eq.~(\ref{ap_peps}), we obtain, at the leading order,
\begin{equation}\label{ap_p0}
p_0^+=p_0^-=p_0,
\end{equation}
and, furthermore, by using Eq.~(\ref{ap_p0}), we also obtain
$$
\partial_xp_0=\Lambda_\delta(M-m)\frac{p_1^--p_1^+}{2},
$$
\begin{equation}\label{ap_j1flux}
\frac{p_1^+-p^-_1}{2}=-\frac{\partial_x p_0}{\Lambda_\delta(M-m)}.
\end{equation}

By taking the sum of Eq. (\ref{ap_peps}), we obtain
$$
\partial_t\left(\frac{p_\varepsilon^++p_\varepsilon^-}{2}\right)
+\partial_x\left(
\frac{p_\varepsilon^+-p_\varepsilon^-}{2\varepsilon}
\right)
+\partial_m\left(
\frac{M-m}{\tilde \tau}\frac{p_\varepsilon^++p_\varepsilon^-}{2}
\right)=0.
$$
Thus, by taking the limit $\varepsilon\rightarrow 0$ in the above equation and using Eq.~(\ref{ap_j1flux}), we obtain
\begin{equation}\label{asymp_largetau}
\partial_t p_0-\partial_x\left(
\frac{\partial_x p_0}{\Lambda_\delta(M(S)-m)}
\right)+\partial_m\left(\frac{M(S)-m}{\tilde \tau}p_0\right)=0.
\end{equation}

The consistency with the Keller-Segel limit can be confirmed when taking the limit $\tilde \tau\rightarrow 0$ in Eq.~(\ref{asymp_largetau}); that is,
at the limit $\tilde \tau\rightarrow 0$, we can obtain from Eq.~(\ref{asymp_largetau})
$$
p_0=\rho_0\delta(m-M(S)),
$$
and
\begin{eqnarray*}
\partial_t \rho_0\delta(m-M(S))-\partial_x\left(
\frac{\partial_x \rho_0\delta(m-M(S))}{\Lambda_\delta(M(S)-m)}
\right)=0,\\
\partial_t \rho_0\delta(m-M(S))-
\partial_{x}
\left(
\frac{\delta(m-M(S))}{\Lambda_\delta(M(S)-m)}\partial_x\rho_0
-\frac{\delta'(m-M(S))}{\Lambda_\delta(M(S)-m)}\rho_0\partial_xM(S)
\right)
=0.
\end{eqnarray*}
Thus, by integrating the above equation with respect to $m$, we obtain Eq.~(\ref{eq_c3_KS}).

\subsection{
Numerical Scheme
}\label{app_scheme}
The one-dimensional space $x\in[-L/2,L/2]$ and the internal state $m\in[-Y,Y]$ are discretized as $x_i=-L/2+i \Delta x$ ($i=0,1,\cdots,I$) and $m_k=-Y+k \Delta m$ ($k=0,1,\cdots,K$), where the mesh intervals are defined as $\Delta x=L/I$ and $\Delta m=2Y/K$.

By integrating Eq.~(\ref{asymp_largetau}) over the unit cell $[x_{i-\frac 12},x_{i+\frac 12}]\times[m_{k-\frac 12},m_{k+\frac 12}]$ and time interval $[t_n,t_{n+1}]$, where $t_n=n\Delta t$, we obtain
\begin{eqnarray*}
p_{i,k}^{n+1}=p_{i,k}^n
&+\frac{1}{\Delta x \Delta m}
\left[
\int_{t_n}^{t_{n+1}}dt \int_{m_{k-\frac 12}}^{m_{k+\frac 12}}dm
\frac{\partial_x p_0}{\Lambda_\delta(M(S)-m)}
\right]_{x_{i-\frac 12}}^{x_{i+\frac 12}}
\\
&-\frac{1}{\Delta x\Delta m}\left [
\int_{t_n}^{t_{n+1}}dt \int_{x_{i-\frac 12}}^{x_{i+\frac 12}}dx
\frac{M(S)-m}{\tilde \tau}p_0
\right]_{m_{k-\frac 12}}^{m_{k+\frac 12}},
\end{eqnarray*}
where $p_{i,k}^n$ is the value of $p_0$ at unit cell $[x_{i-\frac 12},x_{i+\frac 12}]\times[m_{k-\frac 12},m_{k+\frac 12}]$ at time $t_n$.

When we approximate $\partial_x p_0$ in the second term of the R.H. S by the centered difference and apply the upwind scheme for the last term, we obtain
\begin{equation}
p_{i,k}^{n+1}=p_{i,k}^n
+\frac{\Delta t}{\Delta x^2}
\left(
\frac{p^n_{i+1,k}-p^n_{i,k}}{\Lambda_{i+\frac 12,k}}
-\frac{p^n_{i,k}-p^n_{i-1,k}}{\Lambda_{i-\frac 12,k}}
\right)
-\frac{\Delta t}{\tau \Delta m}\left (
\psi_{i,k+\frac 12}-\psi_{i,k-\frac 12}
\right),
\end{equation}
where $\Lambda_{i-\frac 12,k}=\Lambda_\delta(M(S_{i-\frac 12})-m_k)$ and the flux $\psi_{i,k-\frac 12}$ is defined as
%
%
\begin{equation}
    \psi_{i,k-\frac 12}=
    \left(M(S_i)-m_{k-1}\right)^+p_{i,k-1}
    -\left(M(S_i)-m_k\right)^-p_{i,k}.
\end{equation}
Here, we use the notation $u^+=\max\{0,u\}$ and $u^-=\max\{0,-u\}$.

\bibliographystyle{amsxport}
\bibliography{twovel}

\end{document}